\documentclass[11pt, oneside, reqno]{amsart}  	

\usepackage[colorinlistoftodos]{todonotes}
\usepackage[T1]{fontenc}
\usepackage[utf8]{inputenc}

\usepackage{graphicx}

\usepackage{geometry}                		
\usepackage[parfill]{parskip}    		
\usepackage{amssymb}
\usepackage{graphicx}				
\usepackage{amscd, xcolor, changepage,cancel}
\usepackage{faktor}
\usepackage{bbold}
\usepackage{comment}
\usepackage{centernot}

\usepackage[normalem]{ulem}

\usepackage[breaklinks=true]{hyperref}
\hypersetup{
    colorlinks=true,
    linkcolor=blue,
    citecolor=olive,
    filecolor=blue,      
    urlcolor=black,
}

\usepackage{accents}

\newcommand{\doubletilde}[1]{\accentset{\approx}{#1}}

\usepackage{stackengine}
\stackMath

\newcommand\tripletilde[1]{%
  \stackon[-1.5pt]{\stackon[-1.5pt]{\stackon[-1.5pt]{#1}{\sim}}{\sim}}{\sim}%
}

\usepackage[backend=biber,style=alphabetic,sorting=nyt]{biblatex}
\makeatletter
\def\@tocline#1#2#3#4#5#6#7{\relax
  \ifnum #1>\c@tocdepth 
  \else
    \par \addpenalty\@secpenalty\addvspace{-5pt}%
    \begingroup \hyphenpenalty\@M
    \@ifempty{#4}{%
      \@tempdima\csname r@tocindent\number#1\endcsname\relax
    }{%
      \@tempdima#4\relax
    }%
    \parindent\z@ \leftskip#3\relax \advance\leftskip\@tempdima\relax
    \rightskip\@pnumwidth plus4em \parfillskip-\@pnumwidth
    #5\leavevmode\hskip-\@tempdima
      \ifcase #1
       \or\or  \hskip 1em \or \hskip 2em \else \hskip 3em \fi%
      #6\nobreak\relax
    \hfill\hbox to\@pnumwidth{\@tocpagenum{#7}}\par
    \nobreak
    \endgroup
  \fi}
\makeatother

\definecolor{dimgray}{rgb}{0.41, 0.41, 0.41}
\definecolor{afb}{rgb}{0.36, 0.54, 0.66}
\definecolor{bncs}{rgb}{0.0, 0.53, 0.74}
\definecolor{tangerine}{rgb}{0.95, 0.52, 0.0}

\definecolor{bdf}{rgb}{0.19,0.55, 0.91} 
\definecolor{blush}{rgb}{0.87, 0.36, 0.51}

\DeclareMathOperator{\supp}{supp}

\DeclareMathOperator{\Op}{Op}

\DeclareMathOperator{\RL}{RL}
\DeclareMathOperator{\BCH}{BCH}

\newcommand{\beq}{\begin{equation}}
\newcommand{\eeq}{\end{equation}}

\newcommand{\beqq}{\begin{equation*}}
\newcommand{\eeqq}{\end{equation*}}

\newcommand{\mS}{\mathcal{S}}
\newcommand{\mB}{\mathcal{B}}
\newcommand{\mE}{\mathcal{E}}

\newcommand{\mL}{\mathcal{L}}
\newcommand{\mJ}{\mathcal{J}}

\newcommand{\mD}{\mathcal{D}}
\newcommand{\mC}{\mathcal{C}}

\newcommand{\R}{\mathbb{R}}
\newcommand{\C}{\mathbb{C}}
\newcommand{\N}{\mathbb{N}}
\newcommand{\Z}{\mathbb{Z}}

\newcommand{\fg}{\mathfrak{g}}

\newcommand{\la}{\langle}
\newcommand{\ra}{\rangle}

\newcommand{\norm}[1]{\|#1\|}
\newcommand{\RLs}[2]{\|#1\|_{\RL^2_{\lp#2\rp}}}

\newcommand{\ap}{\alpha}
\newcommand{\vp}{\varphi}
\newcommand{\de}{\delta}

\newcommand{\g}{\gamma}

\newcommand{\ep}{\epsilon}

\newcommand{\wt}{\widetilde}
\newcommand{\wtp}{\widetilde{\psi}}

\newcommand{\wte}{\widetilde{\eta}}

\newcommand{\p}{\partial}

\newcommand{\lp}{\left(}
\newcommand{\rp}{\right)}

\newcommand{\lf}{\left|}
\newcommand{\rf}{\right|}

\numberwithin{equation}{section}
\allowdisplaybreaks

\newtheorem{theorem}{Theorem}[section]

\newtheorem{proposition}[theorem]{Proposition}

\newtheorem{lemma}[theorem]{Lemma}

\theoremstyle{definition}
\newtheorem{definition}[theorem]{Definition}

\newtheorem{claim}[theorem]{Claim}
\theoremstyle{remark}
\newtheorem{remark}[theorem]{Remark}

\makeatletter
\@namedef{subjclassname@2020}{%
  \textup{2020} Mathematics Subject Classification}
\makeatother

\title[Inverses of product and flag kernels]{Inverses of product kernels and flag kernels on graded Lie groups}

\author[A. Stokolosa]{Amelia Stokolosa}
\address{Georgia Institute of Technology} \email{stokolosa@gatech.edu}

\begin{document}


\maketitle

\begin{abstract}
Let $T (f) = f*K$, where $K$ is a product kernel or a flag kernel on a direct product of graded Lie groups $G = G_1 \times \cdots \times G_{\nu}$. Suppose $T$ is invertible on $L^2(G)$. We prove that its inverse is given by $T^{-1}(g) = g*L$, where $L$ is a product kernel or a flag kernel accordingly, thus extending a single-parameter inversion theorem by Christ and Geller. 
\end{abstract}


\section{Introduction}

We establish an inversion theorem extending a single-parameter result of Christ and Geller in \cite{CG84} to the multi-parameter setting. Christ and Geller considered operators given by the non-commutative group convolution $Tf = f*K$ on a graded Lie group $G$, where $K$ is a homogeneous kernel with respect to single-parameter dilations. Our inversion theorem applies to a larger class of kernels $K$, defined on a direct product of graded Lie groups $G = G_1 \times \cdots \times G_{\nu}$, which are \textit{almost homogeneous} with respect to \textit{multi-parameter dilations}, namely product kernels and flag kernels (see Definition \ref{def PK} and Definition \ref{def FK}).

R. Fefferman and Stein in \cite{FS82}, and Journ\'e in \cite{Jou85} first introduced product singular integral operators on Euclidean product spaces. Flag singular integral operators appeared later on in the work of Müller, Ricci, and Stein in their study of spectral multipliers on Heisenberg-type groups in \cite{MRS95}. They also obtained the $L^p$ boundedness of operators $Tf = f*K$, where $K$ is a product kernel on the direct product of two stratified Lie groups $G = G_1 \times G_2$ with a biparameter structure. Nagel, Ricci, and Stein investigated the general multi-parameter case while searching for estimates for fundamental solutions of the Kohn-Laplacian $\Box_b$. In particular, they considered operators $Tf = f*K$, where $K$ is a product kernel or a flag kernel on a direct product of homogeneous nilpotent Lie groups $G= G_1 \times \cdots \times G_{\nu}$ (see \cite{NRS01}).  The theory of such operators and their variants quickly developed thereafter and found many applications (see \cite{CF85}, \cite{RS92}, \cite{NS04}, \cite{Yan09}, \cite{Glo10}, and more recently \cite{WL12}, \cite{NRSW12}, \cite{Glo13}, \cite{NRSW18}, \cite{DLOUPW19}).

In the interest of clarity, consider the $2$-parameter setting (we refer the reader to section 2 for a description of the general $\nu$-parameter setting). Let $\fg_1$ be a finite-dimensional graded Lie algebra. By definition, $\fg_1$ decomposes into a direct sum of vector spaces; that is, for some integer $n_1 \in \N$, we have
\begin{align*}
    \fg_1 = \bigoplus_{l=1}^{n_1} V_l^1,
\end{align*}
where $[V_{l_1}^1, V_{l_2}^1] \subseteq V_{l_1 + l_2}^1$ and $V_l^1 = \{0\}$ for $l > n_1$. The exponential map $\exp: \fg_1 \rightarrow G_1$, where $G_1$ is the associated connected, simply connected graded Lie group, is a diffeomorphism\footnote{See Proposition 1.2 in \cite{FS82} for a proof of this statement.}. We henceforth identify $G_1$ with $\R^{q_1}$, where $q_1 = \sum_{l=1}^{n_1} q_l^1$ and $q_l^1= \dim V^1_l$. Notice that $\R^{q_1}$ inherits a non-commutative group multiplication which one can compute explicitly via the Baker-Campbell-Hausdorff formula. With this construction, we define single-parameter non-isotropic dilations on $\R^{q_1}$: for $r_1>0$ and $t_1 = (t_1^1, \ldots, t_{n_1}^1) \in \R^{q_1} = \R^{q_1^1} \times \cdots \times \R^{q^1_{n_1}}$, we define
\vspace{-.1in}
\begin{align*}
    r_1 \cdot t_1 = (r_1t_1^1, r_1^2 t_2^1, \ldots, r_1^{n_1}t_{n_1}^1).
\end{align*}
Let $Q_1 = \sum_{l=1}^{n_1}l \cdot q_l^1$ denote the associated ``homogeneous dimension'' of $G_1$. Similarly, let $\fg_2$ be another finite-dimensional graded Lie algebra with an associated graded Lie group $G_2$ which we identify with $\R^{q_2}$. We thus obtain a direct product of graded Lie groups $G = G_1 \times G_2$ which we identify with $\R^q = \R^{q_1} \times \R^{q_2}$, where $q= q_1 + q_2$. Finally, we define a family of $2$-parameter dilations on $G$ as follows. For $r=(r_1, r_2) \in [0,\infty)^2$, let $r\cdot t = (r_1 \cdot t_1, r_2 \cdot t_2)$.

Product kernels relative to the decomposition $\R^{q_1} \times \R^{q_2}$ are distributions satisfying a growth condition: given a multi-index $(\ap_1, \ap_2) \in \N^{q_1} \times \N^{q_2}$,
\beq
    \lf \p_{t_1}^{\ap_1} \p_{t_2}^{\ap_2} K(t_1, t_2) \rf \leq C_{(\ap_1, \ap_2)} |t_1|_1^{-Q_1 - \deg \ap_1} |t_2|_2^{-Q_2 - \deg \ap_2} ,
\eeq
where $|\cdot|_{\mu}$ is a ``homogeneous norm'' on $\R^{q_{\mu}}$, for $\mu \in \{ 1, 2 \}$ (see an explicit formula for $|\cdot|_{\mu}$ in \eqref{eq hom norm}). In particular, product kernels are smooth away from the ``cross'' $t_1 =0$, $t_2 =0$. They also satisfy a cancellation condition defined recursively (see Definition \ref{def PK}). On the other hand, flag kernels satisfy a growth condition that presents more singularity in the first variable:
\begin{equation}
    |\p_{t_1}^{\ap_1}  \p_{t_2}^{\ap_2} K(t_1, t_2)| \leq C_{(\ap_1, \ap_2)} |t_1|_1^{-Q_1 - \deg \ap_1} \lp |t_1|_1 + |t_2|_2 \rp^{-Q_2 - \deg \ap_2}. 
\end{equation}
Flag kernels are thus smooth away from the coordinate axis $t_1=0$. They also satisfy a cancellation condition defined recursively (see Definition \ref{def FK}). Our main result is as follows:
\begin{theorem}\label{main thm PK}
Let $T$ be a left-invariant singular integral operator given by $T(f) = f*K$, where $K$ is a product kernel (respectively a flag kernel) on a direct product of graded Lie groups $G= G_1 \times \cdots \times G_{\nu}$. If $T$ is invertible as an operator on $L^2(G)$, then its inverse is also of the form $T^{-1}(g) = g*L$, where $L$ is a product kernel (respectively a flag kernel).
\end{theorem}

Most operations and operators on $G_1 \times G_2$ do not commute. For example, group multiplication and group convolution are both non-commutative. Nonetheless, right-invariant operators and left-invariant operators commute by associativity of convolution. As such, to prove regularity properties of the inverse $T^{-1}$, we introduce \textit{right-invariant} differential operators $(I+\mL_{\mu})^{s_{\mu}}$ on each factor space $G_{\mu}$ with which the \textit{left-invariant} operator $T$ commutes:
\beq
    (I+\mL_{\mu})^{s_{\mu}} T  = T (I+\mL_{\mu})^{s_{\mu}},
\eeq
for $\mu =1, 2$. 
To construct the central ideas in our proof, we extend a single-parameter \textit{a priori estimate} by Christ and Geller (see Lemma 5.3 in \cite{CG84}) to the multi-parameter setting. The key idea in our proof are the \textit{a priori estimates} in Proposition \ref{CG lemma PK} and Proposition \ref{CG lemma FK} which apply to a larger class of not necessarily homogeneous multi-parameter singular integrals, namely operators given by group convolution on the right with a product kernel or a flag kernel.

\begin{remark}
K{\c e}pa obtained a related inversion theorem for flag kernels on the Heisenberg group in \cite{K16} using representation theory. We use tools from PDEs instead of representation theory to construct an \textit{a priori estimate} in Proposition \ref{CG lemma FK} for flag kernels defined on a direct product of graded Lie groups. Other notable single-parameter inversion theorems include the foundational work by Calderón and Zygmund in \cite{CZ56} on Euclidean spaces, and the works on inverses of single-parameter singular kernels by \cite{Ch88}, \cite{Ch88b}, \cite{CGGP92}, \cite{Wei08}, and \cite{Glo17}. 
\end{remark}

\vspace{.2in}

\section{Background and Notation}

For every $\mu \in \{ 1, \ldots, \nu \}$, let $\fg_{\mu}$ be a finite-dimensional graded Lie algebra. By definition, $\fg_{\mu}$ decomposes into a direct sum of vector spaces
\begin{align*}
    \fg_{\mu} = \bigoplus_{l=1}^{n_{\mu}} V_l^{\mu}, 
\end{align*}
where $[V_{l_1}^{\mu}, V_{l_2}^{\mu}] \subseteq V_{l_1 + l_2}^{\mu}$ and $V_l^{\mu} = \{0\}$, for $l > n_{\mu}$. For every $l \in \{ 1, \ldots, n_{\mu} \}$, let $\{ X_{k_l }^{\mu}\}_{k_l =1 }^{q_l^{\mu}}$ be a basis of left-invariant vector fields for $V_l^{\mu}$ so that $q_l^{\mu} = \dim V_l^{\mu}$. In addition, let $q_{\mu} = \sum_{l=1}^{n_{\mu}} q_l^{\mu}$ and let $\{X_1^{\mu}, \ldots, X_{q_{\mu}}^{\mu}\}$ be an enumeration of these basis vector fields, thereby forming a basis for $\fg_{\mu}$. The exponential map $\exp: \fg_{\mu} \rightarrow G_{\mu}$, where $G_{\mu}$ is the associated connected, simply connected graded Lie group, is a diffeomorphism. We thus obtain global coordinates $\R^{q_{\mu}} \rightarrow G_{\mu}$: 
\begin{align*}
    (t_1^{\mu}, \ldots, t_{q_{\mu}}^{\mu}) \mapsto \exp(t_1^{\mu} X_1^{\mu} + \ldots + t^{\mu}_{q_{\mu}} X_{q_{\mu}}^{\mu} ). 
\end{align*}
Given $x_{\mu} = (x_1^{\mu}, \ldots, x_{q_{\mu}}^{\mu})$ and $ y_{\mu}= (y_1^{\mu}, \ldots, y_{q_{\mu}}^{\mu}) \in \R^{q_{\mu}}$, one can compute the group multiplication $x_{\mu} \cdot y_{\mu}$ which is given by the coefficients of the basis vectors after applying the Baker-Campbell-Hausdorff formula\footnote{The Baker-Campbell-Hausdorff formula is given by $
		\BCH(x^{\mu} \cdot X^{\mu}, y^{\mu} \cdot X^{\mu}) = \exp \Big( x^{\mu} \cdot X^{\mu} +y^{\mu} \cdot X^{\mu} + \frac{1}{2} [x^{\mu} \cdot X^{\mu}, x^{\mu} \cdot X^{\mu}] + \frac{1}{12} [x^{\mu} \cdot X^{\mu},[x^{\mu} \cdot X^{\mu},y^{\mu} \cdot X^{\mu}]] $ $ + \frac{1}{12} [y^{\mu} \cdot X^{\mu},[y^{\mu} \cdot X^{\mu}, x^{\mu} \cdot X^{\mu}]]  + \textrm{higher order commutators} \Big).$ It is an identity and not just a formal series as all commutators of high enough order vanish. See equation  (2.73) in \cite{Str14} for reference. }:
\begin{align*}
    \BCH(x_1^{\mu} X_1^{\mu} + \ldots + x^{\mu}_{q_{\mu}} X_{q_{\mu}}^{\mu}, y_1^{\mu} X_1^{\mu} + \ldots + y^{\mu}_{q_{\mu}} X_{q_{\mu}}^{\mu}).
\end{align*}  
We henceforth identify $G_{\mu}$ with $\R^{q_{\mu}} = \R^{q_1^{\mu}} \times \cdots \times \R^{q_{n_{\mu}}^{\mu}}$ and obtain a family of automorphisms, called \textit{single-parameter dilations}, on $\R^{q_{\mu}}$: for $r_{\mu}>0$, let
\begin{align}\label{eq SP dilations}
    r_{\mu} \cdot t_{\mu} = (r_{\mu} t^{\mu}_1, r_{\mu}^2 t^{\mu}_2, \ldots, r_{\mu}^{n_{\mu}} t^{\mu}_{n_{\mu}}). 
\end{align}
Let $Q_{\mu} = \sum_{l=1}^{n_{\mu}}l \cdot q_l^{\mu}$ be the associated ``homogeneous dimension'' of $\R^{q_{\mu}}$. 

\begin{definition}
A \textit{homogeneous norm} on $\R^{q_{\mu}}$ is a continuous function $|\cdot|_{\mu}: \R^{q_{\mu}} \rightarrow [0,\infty)$ that is smooth away from $0$ with $|t_{\mu}|_{\mu} =0$ $\Leftrightarrow \ t_{\mu} =0$ and $|r_{\mu} \cdot t_{\mu}|_{\mu} = r_{\mu}|t_{\mu}|_{\mu}$, for $r_{\mu}>0$.
\end{definition}

Any two such homogeneous norms on $\R^{q_{\mu}}$ are equivalent. Given $X = \sum_{l=1}^{n_{\mu}} \sum_{k_l =1}^{ q_l^{\mu} } t_{l, k_l}^{\mu} X_{k_l}^{\mu}$, we thus define
\beq\label{eq hom norm}
    |t_{\mu}|_{\mu}:= \Bigg( \sum_{l=1}^{n_{\mu}} \sum_{k_l =1}^{q_l^{\mu} } |t_{l, k_l}^{\mu}|^{2(n_{\mu}!)/ l} \Bigg)^{1/(2 (n_{\mu}!))}.
\eeq
Let $\{X_1^{\mu}, \ldots, X_{q_{\mu}}^{\mu}\}$ and $\{Y_1^{\mu}, \ldots, Y_{q_{\mu}}^{\mu}\}$ be spanning sets of left- and right-invariant vector fields on $G_{\mu}$ s.t. at the identity, $X_j^{\mu} = Y_j^{\mu} = \frac{\p}{\p x_j^{\mu}}$. Note that $X_j^{\mu}$ and $Y_j^{\mu}$ are both homogeneous\footnote{That is, for all $r_{\mu} >0$, $D(f(r_{\mu} \cdot t_{\mu})) = r_{\mu}^l (Df)(r_{\mu} \cdot t_{\mu})$, where $D=$ $X_j^{\mu}$, $Y_j^{\mu}$. } of degree $l$, provided $x_j^{\mu} \in \R^{q_l^{\mu}}$.

For $r \in [0,\infty)^{\nu}$, we define multi-parameter dilations using the single-parameter dilations defined in \eqref{eq SP dilations} on each factor space: 
\begin{align}\label{eq MP dilations}
    r \cdot t = (r_1 \cdot t_1, \ldots, r_{\nu} \cdot t_{\nu}). 
\end{align}
In addition, let $rX$ denote the following ordered list of vector fields with appropriate dilations:
\begin{equation*}
    rX= r_1  X^1, \ldots, r_{\nu}  X^{\nu}
    = r_1^{d_1^1}X_1^1, \ldots , r_1^{d_{q_1}^1} X_{q_1}^1, \ldots, r_{\nu}^{d_1^{\nu}} X_1^{\nu}, \ldots, r_{\nu}^{ d_{q_{\nu}}^{\nu}}X_{q_{\nu}}^{\nu},
\end{equation*}
where $d_j^{\mu}=l$, provided $X_j^{\mu} \in V_l^{\mu}$ where $l\in \{1, \ldots, n_{\mu}\}$, for every $j \in \{ 1, \ldots, q_{\mu} \}$ and $\mu \in \{ 1, \ldots, \nu \} $. 

For every multi-index $\ap_{\mu} \in \N^{q_{\mu}} = \N^{q_1^{\mu}} \times \cdots \times \N^{q_{n_{\mu}}^{\mu}}$, let $\deg \ap_{\mu} = \sum_{l=1}^{n_{\mu}} l \norm{\ap_l^{\mu}}_{l^1}$ denote its homogeneous degree and $|\ap_{\mu}|= \sum_{l=1}^{n_{\mu}} \norm{\ap_l^{\mu}}_{l^1}$, its isotropic degree. In addition, for every multi-index $\ap = (\ap_1, \ldots, \ap_{\nu}) \in \N^{q_1} \times \cdots \times \N^{q_{\nu}}$, let $|\ap|= (|\ap_1|, \ldots, |\ap_{\nu}|)$ and $\deg \ap = \lp \deg \ap_{1}, \ldots, \deg \ap_{\nu} \rp$.

\vspace{.2in}

\section{Inversion Theorem for Product Kernels}

\begin{definition}\label{def PK}
A \textit{product kernel} $K$ on $\R^q$, relative to the decomposition $\R^q = \R^{q_1} \times \cdots  \times \R^{q_{\nu}}$, is a distribution which coincides with a  $C^{\infty}$ function away from the coordinate subspaces $t_{1} =0, \ldots, t_{\nu} =0$ and which satisfies the following two conditions:

\begin{enumerate}
    \item Growth condition - For each multi-index $\ap = (\ap_1, \ldots, \ap_{\nu}) \in \N^{q_1} \times \cdots  \times \N^{q_{\nu}}$, there exists a constant $C_{\ap}$ such that, away from the coordinate subspaces $t_1=0$, $\ldots$, $t_{\nu} =0$, 
\begin{equation}\label{GC PK}
|\p_{t_1}^{\ap_1} \cdots \p_{t_{\nu}}^{\ap_{\nu}} K(t)| \leq C_{\ap} |t_1|_1^{-Q_1 - \deg \ap_1} \cdots |t_{\nu}|_{\nu}^{-Q_{\nu} - \deg \ap_{\nu}}.
\end{equation}
For every multi-index $\ap$, we take the least $C_{\ap}$ to define a seminorm. 

    \item Cancellation condition - This condition is defined recursively. 

    \begin{itemize}
        \item For $\nu =1$, given a bounded set\footnote{As a corollary of Proposition 14.6 p.139 in \cite{Tr67}, a set $\mB \subseteq C^{\infty}_0(\R^n)$ is bounded if the following two conditions hold: 

(1) there exists a compact set $K \Subset \R^{n}$ s.t. for all $f \in \mB$, $\supp f \subseteq K$;

(2) for every multi-index $\ap \in \N^n$, $\sup_{x\in \R^{n};  f\in \mB} | \p^{\ap} f(x)|<\infty$,

\noindent where $C^{\infty}_0(\R^n)$ denotes the set of compactly supported smooth functions.} $\mB \subseteq C^{\infty}_0(\R^q)$, 
\begin{equation}\label{CC1 PK}
\sup_{\phi \in \mB; R>0} \lf \int K(t) \phi(R\cdot t) dt \rf <\infty. 
\end{equation}

        \item For $\nu >1$, given $1 \leq \mu \leq \nu$, a bounded set $\mB_{\mu} \subseteq C^{\infty}_0(\R^{q_{\mu}})$, $\phi_{\mu} \in \mB_{\mu}$, and $R_{\mu}>0$, the distribution $K_{\phi_{\mu},R_{\mu}}$ defined by
\begin{equation}\label{CC2 PK}
    K_{\phi_{\mu}, R_{\mu}} (\ldots, t_{\mu -1},  t_{\mu+1}, \ldots) := \int K(t)\phi_{\mu}(R_{\mu}\cdot  t_{\mu})dt_{\mu}
\end{equation}
is a product kernel on the $(\nu-1)$-factor space $ \cdots \times \R^{q_{\mu -1}} \times \R^{q_{\mu+1}} \times \cdots  $, where the bounds are independent of the choice of $\phi_{\mu}$ and $R_{\mu}$. 
    \end{itemize}

\noindent For the base case $\nu=0$, we define the space of product kernels to be $\C$ with its usual topology. For $\nu \geq 1$, given a seminorm $|\cdot|$ on the space of $(\nu-1)$-factor product kernels, we define a seminorm on the $\nu$-factor product kernels by 
\begin{equation}\label{seminorm2 PK}
    |K|:= \sup_{\phi_{\mu} \in \mB_{\mu}; R_{\mu}>0} |K_{\phi_{\mu}, R_{\mu}}|,
\end{equation}
which we assume to be finite.
\end{enumerate}

\end{definition}

\begin{remark}
\cite{FS82} and \cite{Jou85} introduced product singular integral operators on Euclidean spaces. \cite{MRS95} later on defined product kernels $K$ on the direct product of homogeneous groups $G = G_1 \times G_2$ and proved the $L^p$ boundedness\footnote{See Theorem 4.4 p.221 in \cite{MRS95}.} of the associated left-invariant operator $Tf = f*K$, for $1<p< \infty$. Finally, \cite{NRS01} investigated product kernels on the direct product of homogeneous groups while studying solutions of the Kohn-Laplacian $\Box_b$. 
\end{remark}

\begin{definition}
    A \textit{Calderón-Zygmund kernel} is defined to be a single-parameter product kernel. 
\end{definition}

\begin{definition}\label{bump functions}
A \textit{bounded set of bump functions} on $\R^{q_{\mu}}$ is a set of triples $\{(\phi_{\mu}, z_{\mu}, r_{\mu})\} \subseteq C^{\infty}_0(\R^{q_{\mu}}) \times \R^{q_{\mu}} \times (0,\infty)$ s.t. $\phi_{\mu}(t_{\mu}):= r_{\mu}^{-Q_{\mu}} \psi_{\mu}(r_{\mu}^{-1}\cdot (z_{\mu}^{-1}t_{\mu}))$ where $\{\psi_{\mu}\} \subseteq C^{\infty}_0(B^{\mu}(0,1))$ is a bounded set\footnote{$B^{\mu}(0,1) \subseteq \R^{q_{\mu}}$ denotes the unit ball centered at the identity $0$ in $\R^{q_{\mu}}$ with respect to the homogeneous norm $|\cdot|_{\mu}$.}.
\end{definition}

\begin{definition}\label{def PSIO} An operator $S: C^{\infty}_0(\R^q) \rightarrow C^{\infty}(\R^q)$ is a \textit{product singular integral operator of order $s = (s_1, \ldots , s_{\nu})$} $\in (-Q_1, \infty) \times \ldots \times (-Q_{\nu}, \infty)$ if it satisfies the following conditions:

\begin{enumerate}
    \item Growth Condition - For all multi-indices $\ap, \beta$, there exists $C_{\ap, \beta}$ s.t.
\begin{equation}\label{GC PSIO}
    |X^{\ap}_x X^{\beta}_y S(x,y)| \leq C_{\ap,\beta} \prod_{\mu=1}^{\nu}|y_{\mu}^{-1}x_{\mu}|_{\mu}^{-s_{\mu} - Q_{\mu} -\deg \ap_{\mu} - \deg \beta_{\mu} },
\end{equation}
where $S(x,y)$ denotes the Schwartz kernel of the operator $S$. The least possible $C_{\ap, \beta}$ defines a seminorm.

    \item Cancellation Condition - For every $1 \leq \mu \leq \nu$, and for all bounded sets of bump functions $\{(\phi_{\mu}, z_{\mu}, r_{\mu})\} \subseteq C^{\infty}_0(\R^{q_{\mu}}) \times \R^{q_{\mu}} \times (0, \infty)$, we define a map $x_{\mu} \mapsto S^{\phi_{\mu}, x_{\mu}}$ s.t. 
\begin{equation*}
    \int_{\R^{q_{\mu}}} \la S^{\phi_{\mu},x_{\mu}} (\otimes_{\g \neq \mu} \phi_{\g}), \otimes_{\g \neq \mu} \psi_{\g} \ra \psi_{\mu}(x_{\mu}) dx_{\mu}:=\la S(\phi_1 \otimes \ldots \otimes \phi_{\nu}), \psi_1 \otimes \ldots \otimes \psi_{\nu} \ra.
\end{equation*}
In addition, we assume that for every $\ap$, the operator 
\begin{align*}
    r_{\mu}^{Q_{\mu}+s_{\mu}} (r_{\mu}X_{x_{\mu}}^{\mu})^{\ap} S^{\phi_{\mu},x_{\mu}}
\end{align*}
is a product singular integral operator of order $(\ldots , s_{\mu-1}, s_{\mu+1}, \ldots)$ on the $(\nu-1)$-factor space $\cdots \times \R^{q_{\mu-1}}\times \R^{q_{\mu+1}} \times \cdots $. 

\end{enumerate}
Finally, for every continuous seminorm $|\cdot|$ on the space of product singular integral operators defined on the $(\nu-1)$-factor space $\cdots \times \R^{q_{\mu-1}}\times \R^{q_{\mu+1}} \times \cdots $, for every multi-index $\ap$, and every bounded set of bump functions $\mC_{\mu}$ on $\R^{q_{\mu}}$, we define a seminorm $|\cdot|_{\ap, \mC_{\mu}}$ on product singular integral operators $S$ on $\R^{q_1} \times \cdots \times \R^{q_{\nu}}$ by
\begin{equation}
    \lf S \rf_{\ap, \mC_{\mu}}:= \sup_{\substack{(\phi_{\mu}, z_{\mu}, r_{\mu}) \in \mC_{\mu}; \\ x_{\mu}\in \R^{q_{\mu}}}} \lf r_{\mu}^{Q_{\mu}+s_{\mu}} (r_{\mu}X_{\mu}^{\mu})^{\ap}T^{\phi_{\mu}, x_{\mu}} \rf
\end{equation}
which we assume to be finite. We do the same for the transpose of $S$ in the $\mu$ variable, where we define $z_{\mu} \mapsto S^{z_{\mu}, \psi_{\mu}}$ reversing the roles of $x_{\mu}, z_{\mu}$ and $\phi_{\mu}, \psi_{\mu}$. 

\end{definition}

\begin{remark}
\cite{NS04} introduced multi-parameter product singular integral operators of order $0$ in the sub-Riemannian setting. \cite{Str14} later constructed product singular integral operators of various nonzero order $(s_1, \ldots, s_{\nu})$. 
\end{remark}

\vspace{-.05in}

We will make use of an equivalent definition for product singular integral operators of order $(s_1, \ldots, s_{\nu})$ by Street. To do so, we first introduce the building blocks of such operators in the next definition which adapts Definition 4.1.11 p.228 and Definition 5.1.8 p.270 in \cite{Str14} to our setting.

\begin{definition}\label{def BSEOs X} Let $\Omega := \Omega_1 \times \cdots \times \Omega_{\nu} \Subset \R^{q_1} \times \cdots \times \R^{q_{\nu}}$ be a relatively compact open subset. The \textit{set of bounded sets of elementary operators} $\mathcal{G}$ on $\Omega$ is defined to be the largest set of subsets of $C^{\infty}_0(\Omega \times \Omega) \times (0, 1]^{\nu}$ s.t. for all $\mE \in \mathcal{G}$,
\begin{itemize}
    \item $\forall \ap, \beta, m$, $\exists C$ s.t. $\forall (E_j, 2^{-j}) \in \mE$,
\begin{equation}\label{BSPEOs X}
    |(2^{-j} X_x)^{\ap} (2^{-j}   X_y)^{\beta} E_j(x, y) |   \leq C \prod_{\mu=1}^{\nu}\frac{(1+2^{j_{\mu}}|y_{\mu}^{-1}x_{\mu}|_{\mu})^{-m_{\mu}}}{(2^{-j_{\mu}}+|y_{\mu}^{-1}x_{\mu}|_{\mu})^{Q_{\mu}}},
\end{equation}
where $E_j(x, y)$ denotes the Schwartz kernel of the operator $E_j$.
    \item Let $e=(1, \ldots, 1) \in \N^{\nu}$. $\forall (E_j, 2^{-j}) \in \mE$, we have
\begin{equation}\label{BSEOs X}
    E_j = \sum_{\substack{|\ap|, |\beta| \leq e}} 2^{-\lp 2e- |\ap|-  |\beta| \rp \cdot j} (2^{-j}X)^{\ap} E_{j, \ap, \beta} (2^{-j} X)^{\beta},
\end{equation}
where $\{(E_{j, \ap, \beta}, 2^{-j}); (E_j, 2^{-j}) \in \mE \} \in \mathcal{G}$. 
\end{itemize}
We call elements $\mE \in \mathcal{G}$ \textit{bounded sets of elementary operators} on $\Omega$.    

\end{definition}

\begin{definition}\label{def 2j el op}
    We say $E \in C^{\infty}_0(\Omega \times \Omega)$ is a \textit{$2^{-j}$ elementary operator} if $\{(E, 2^{-j}); j \in \Z_{\geq 0}^{\nu}\}$ is a bounded set of elementary operators. 
\end{definition}

\cite{Str14} presents four equivalent definitions for product singular integral operators in a more general local setting in Theorem 5.1.12 p.271. We record two of the four equivalent definitions in our local ``product setting'' below. 

\begin{theorem}\label{PSIO theorem BSEOs} Let $\Omega \Subset \R^q$ be a relatively compact open subset. Fix $s \in (-Q_1, \infty) \times \cdots \times (-Q_{\nu}, \infty)$. For $S: C^{\infty}(\Omega) \rightarrow C_0^{\infty}(\Omega)$, the following are equivalent:
\begin{itemize}
    \item $S: C^{\infty}(\Omega) \rightarrow C_0^{\infty}(\Omega)$ is a product singular integral operator of order $s$ as in Definition \ref{def PSIO}. 

    \item  There exists a bounded set of elementary operators $\{(E_j, 2^{-j});  j \in \Z^{\nu}_{\geq 0}\}$ s.t. 
\beq\label{eq def 2 PSIO}
    S = \sum_{j \in \Z^{\nu}_{\geq 0}} 2^{j \cdot s} E_j,
\eeq
where the sum converges in the topology of bounded convergence as operators $C^{\infty}(\Omega) \rightarrow C_0^{\infty}(\Omega)$\footnote{For every continuous seminorm $|\cdot|$ on $C_0^{\infty}(\Omega)$ and for every bounded set $\mB \subseteq C_0^{\infty}(\Omega)$, we define a semi-norm $|\cdot|'$ on the space of continuous linear maps $T:C^{\infty}(\Omega)\rightarrow C_0^{\infty}(\Omega)$ by $|T|' = \sup_{f\in \mB} |Tf|$. The coarsest topology according to which the above semi-norms are continuous is called the ``topology of bounded convergence''.}. 
\end{itemize}

\end{theorem}

\begin{definition}
For $S \subseteq \{1, \ldots, \nu\}$, we define the space $\mS_0^{S}$ as follows:
\beqq
    \mS_0^{S}:= \Big\{ f \in \mS;  \forall \mu \in S, \int f(t)  t_{\mu}^{\ap} \  dt_{\mu} =0, \ \forall \ap \in \N^{q_{\mu}} \Big\}.
\eeqq
\end{definition}

\begin{remark}\label{remark pull out derivatives}
For $\zeta \in \mS_0^{S}$, we can ``pull out derivatives'' in $t_{\mu}$, provided $\mu \in S$. That is, $\zeta = \sum_{|\ap_{\mu}|=1} \p_{t_{\mu}}^{\ap_{\mu}}\zeta_{\ap_{\mu}},$ where $\{ \zeta_{\ap_{\mu}}; |\ap_{\mu}|=1\} \subseteq \mS_0^{S}$ is bounded\footnote{One can verify this statement via a straightforward adaptation of Lemma 1.1.16 p.11 in \cite{Str14} to the setting of graded Lie groups.}.
\end{remark}


To avoid notational headaches from juggling numerous indices and to highlight the key ideas in the proof, we focus on the 2-parameter case. The general $\nu$-parameter case follows from a few straightforward modifications.

\subsection{A Multi-parameter A Priori Estimate} \ 

To prove the desired regularity properties of the inverse operator $T^{-1}$, we introduce \textit{right-invariant} non-isotropic Sobolev spaces. Given $\mu  = 1, 2$, we introduce the following homogeneous, nonnegative, and essentially self-adjoint operators on $\R^{q_{\mu}}$: 
\begin{align*}
    \mL_{\mu} := \sum_{j=1}^{q_{\mu}} \lp Y_j^{\mu} \rp^{\frac{4 n_{\mu}! }{d^{\mu}_j}}.
\end{align*}
We define an analytic family of operators $\big\{\mJ^{\mu}_{s}\big\}_{s \in \C}$ on each factor space $\R^{q_{\mu}}$ (see Proposition 5.1 in \cite{CG84} which adapts the constructions in \cite{Fol75} to the graded Lie group setting) so that $\mJ^{\mu}_0= I$, $\mJ^{\mu}_1 = I + \mL_{\mu}$, $\mJ_{s}^{\mu} \circ \mJ_{t}^{\mu} = \mJ_{s+t}^{\mu}$ and $\mJ_{s}^{\mu}: \mS(\R^{q_{\mu}}) \rightarrow \mS(\R^{q_{\mu}})$. In turn, with these right-invariant operators, we define multi-parameter, non-isotropic Sobolev norms on $\R^q$. The operators above commute so we write $\mJ_{(s_1, s_2)} := \mJ_{s_1}^1 \circ \mJ_{s_2}^{2}$.

\begin{definition}
Given $s = (s_1, s_2) \in \R^{2}$, let $\RL^2_{(s_1, s_2)} (\R^q)$ be the completion of $C^{\infty}_0(\R^{q})$ under the norm\footnote{We label these Sobolev norms with a capital ``R'' to highlight their main characteristic: they are defined by right-invariant differential operators that commute with the left-invariant singular integral operator $T$.}
\begin{equation}\label{right sob}
    \RLs{f}{s_1, s_2} := \|\mJ_{(s_1, s_2)} f\|_{L^2}.
\end{equation}
\end{definition}   

\begin{remark}
The single-parameter non-isotropic Sobolev spaces on nilpotent Lie groups were first introduced by Folland, Rothschild, and Stein in \cite{FS74} and \cite{RS76}. 
\end{remark}

Here is the key multi-parameter a priori estimate for product kernels: 
\begin{proposition}\label{CG lemma PK}
    There exist $\ep_{\mu}>0$ such that for all $\psi_{\mu} \prec\footnote{Henceforth, the notation $\phi \prec \gamma$ will mean that $\phi \gamma = \phi$.} \eta_{\mu} \in C^{\infty}_0(\R^{q_{\mu}})$, $l_{\mu} \in \N$, and $f \in C^{\infty}_0(\R^q)$, where $\mu =1, 2$,
\begin{equation*}
\begin{split}
    \RLs{\psi_1 \otimes \psi_2  f}{l_1\ep_1, l_2\ep_2} \lesssim &\RLs{\psi_1 \otimes \psi_2 Tf}{l_1\ep_1, l_2 \ep_2} + \RLs{\eta_1 Tf}{l_1\ep_1, l_2 \ep_2 - \ep_2} \\
    &+ \RLs{\eta_2 Tf}{l_1 \ep_1-\ep_1, l_2\ep_2}  + \RLs{f}{l_1\ep_1-\ep_1,l_2 \ep_2-\ep_2},
    \end{split}
\end{equation*}
where the implicit constant depends on the test functions $\psi_{\mu}, \eta_{\mu}$ and on the operators $T$ and $T^{-1}$ in an admissible\footnote{The constant depends on the seminorms of the original product kernel $|K|$ (see \eqref{GC PK}, \eqref{seminorm2 PK}), and on the two operator norms $\norm{T}_{\mB(L^2)}$ and $\norm{T^{-1}}_{ \mB(L^2) }$. } way. 
\end{proposition}

We first catalog two results by Nagel and Stein, and Street which we use in the proof of Proposition \ref{CG lemma PK}. 

\begin{theorem}[{\cite[Theorem 4.1.2]{NS04}}]\label{PSIO L2 bdd}
Product singular integral operators $T: C^{\infty}_0(\R^{q_1} \times \R^{q_2}) \rightarrow C^{\infty}(\R^{q_1} \times \R^{q_2})$ of order $(0, 0)$ are bounded on $L^p$, for $1<p<\infty$. 
\end{theorem}

The following theorem, adapted to our graded Lie group setting, says that product singular integral operators on a fixed open and relatively compact subset $\Omega = \Omega_1 \times \Omega_2 \Subset \R^q$ form a \textit{filtered algebra}. 

\begin{theorem}[{\cite[Corollary 5.1.13]{Str14}}]\label{PSIO algebra}
    If $T, S: C^{\infty}(\Omega) \rightarrow C_0^{\infty} (\Omega)$ are product singular integral operators of order $t$ and $s$ respectively, then $T\circ S$ is a product singular integral operator of order $t+s$, for $t, s \in \R^2$. 
\end{theorem}

To obtain the key a priori estimate in Proposition \ref{CG lemma PK}, we first establish the following commutator estimates. 
\begin{lemma}\label{comm SP}
For all $s_1, s_2>0$, $\psi_1 \in C^{\infty}_0(\R^{q_1})$, and $f \in C^{\infty}_0(\R^q)$, there exists $\ep_1>0$ s.t.
\begin{equation}\label{comm J SP}
    \RLs{[\psi_1, \mJ_{(s_1,0)}] f}{\ep_1,s_2} \leq C(s_1, \psi_1) \RLs{f}{s_1,s_2}.
\end{equation}
In addition, given $\phi_1, \psi_1 \in C^{\infty}_0(\R^{q_1})$, for any $f \in C^{\infty}_0(\R^q)$, we have
\begin{equation}\label{comm T SP}
            \norm{ \phi_1 \mJ_{(\ep_1,0)} [T,\psi_1] f}_{L^2} \lesssim C(\phi_1, \psi_1, T, T^{-1}) \norm{f}_{L^2},
\end{equation}
where the implicit constant depends on $T$ and $T^{-1}$ in an admissible way\footnote{See the definition of an admissible constant in Proposition \ref{CG lemma PK}. }. 
\end{lemma}

\begin{remark}
By symmetry, we obtain analogous estimates for $[\psi_2, \mJ_{(0, s_2)}]$ and $\phi_2 \mJ_{(0, \ep_2)} [T,\psi_2]$. 
\end{remark}

We in turn need to prove the following technical lemma which is used in the proof Lemma \ref{comm SP}. 

\begin{lemma}\label{lemma T psi1 pk}
Let $\eta, \eta' \in C^{\infty}_0(\Omega)$. There exists a bounded set of elementary operators \\ $\big\{(E_j, 2^{-j});  j \in \Z^2_{\geq 0}\big\}$ s.t. 
\begin{equation}\label{eq [T,p] decomp}
    \eta [T, \psi_1]\eta' = \sum_{\substack{(j_1, j_2)\in \Z^2_{\geq 0}}} 2^{-j_1}E_{(j_1, j_2)},
\end{equation}
where the sum converges in the topology of bounded convergence\footnote{See Lemma 5.3.2 p.293 in \cite{Str14} for a proof of the convergence. } as operators $C^{\infty}(\Omega) \rightarrow C^{\infty}_0(\Omega)$. 
\end{lemma}

\noindent\textit{Notation - }Let $\Op(g)f := f*g$. In addition, given $j= (j_1, j_2) \in \Z^2$, let $f^{(2^j)}(t_1, t_2) := 2^{j_1 Q_1} 2^{j_2 Q_2} f(2^{j_1} \cdot t_1, 2^{j_2} \cdot t_2)$.

\begin{proof}[Proof of Lemma \ref{lemma T psi1 pk}]
    Consider the Littlewood-Paley decomposition of $T$:
\beq\label{decomp PK}
    T = \sum_{j \in \Z^2} D_j:=  \sum_{j \in \Z^2} \Op\Big(\zeta_j^{(2^j)}\Big),
\eeq
where $\big\{\zeta_j;  j \in \Z^2\big\} \subseteq \mS_0^{\{1, 2\}}$ is bounded\footnote{See Corollary 5.2.16 p.289 in \cite{Str14} for a precise formulation and Theorem 2.2.1 in \cite{NRS01} for an analogous decomposition.}. We thus obtain a decomposition:
\begin{equation}\label{four sums PK}
\begin{split}
    \sum_{j\in \Z^2} \eta [D_j, \psi_1] \eta' 
    =&\sum_{\substack{ j_1, j_2 \leq 0}} \eta [D_j, \psi_1] \eta' +\sum_{\substack{ j_1 \leq 0 < j_2 }} \eta [D_j, \psi_1] \eta' \\
    &+\sum_{\substack{j_2 \leq 0<j_1}} \eta [D_j, \psi_1] \eta' + \sum_{\substack{ j_1, j_2 > 0}} \eta [D_j, \psi_1] \eta' .
\end{split}
\end{equation}
We begin by showing that the first sum on the right-hand side of the equal sign in \eqref{four sums PK} converges to a $2^{0}$ elementary operator which we denote $E_0$. To verify that the associated Schwartz kernel $E_0(x, y)$ satisfies the first condition \eqref{BSPEOs X} in Definition \ref{def BSEOs X}, with $|\ap| = |\beta| = 0$, by the triangle inequality, we have 
\begin{align*}
    |E_0(x, y)| &\lesssim \sum_{\substack{ j_1 , j_2 \leq 0}} 2^{j_1Q_1} 2^{j_2 Q_2}  \lf \zeta_j(2^{j_1} \cdot y_1^{-1}x_1, 2^{j_2}\cdot y_2^{-1}x_2) \rf,
\end{align*}
where $\{\zeta_j; j_1, j_2 \leq 0\} \subseteq \mS$ is bounded. Hence, for all $m_{\mu} >0$, where $\mu =1, 2$,
\begin{align*}
    |E_0(x,y)| \lesssim \sum_{\substack{j_1, j_2\leq 0}} \prod_{\mu =1 ,2} 2^{j_{\mu}Q_{\mu}}  (1+2^{j_{\mu}}|y_{\mu}^{-1}x_{\mu}|_{\mu})^{-m_{\mu}}. 
\end{align*}
Recall $x, y \in \Omega$, a bounded set, and $j_{\mu} \leq 0$. We have
\begin{align*}
    |E_0(x,y)|\lesssim \prod_{\mu =1 ,2}   (1+|y_{\mu}^{-1}x_{\mu}|_{\mu})^{-m_{\mu}}.
\end{align*}
Condition \eqref{BSPEOs X} for $|\ap|, |\beta| \neq 0$ follows directly from an application of the Leibniz rule. In addition, $(E_0, 2^{0})$ immediately satisfies the second condition \eqref{BSEOs X} for elementary operators by letting $E_{0, \ap, \beta} \equiv 0$ whenever $|\ap|+|\beta|>0$. 

In the next step, we show that the second term on the right-hand side of \eqref{four sums PK} is a sum of $2^{-(0, j_2)}$ elementary operators which we denote $E_{(0, j_2)}$; that is, let 
\beqq
    \sum_{0< j_2} E_{(0,j_2)}:= \sum_{0< j_2} \Bigg( \sum_{j_1 \leq 0} \eta [D_j, \psi_1]\eta' \Bigg). 
\eeqq
We first verify that $\big\{(E_{(0,j_2)}, 2^{-(0,j_2)}); j_2>0\big\}$ satisfies \eqref{BSPEOs X}. By the Leibniz rule again, it suffices to consider the case where $|\ap|= |\beta|=0$. By the triangle inequality, 
\begin{align*}
    |E_{(0,j_2)} (x,y)| \leq \sum_{j_1 \leq 0} |\eta(x) (\psi_1(y_1)-\psi_1(x_1))\zeta_j^{(2^j)}(y^{-1}x) \eta'(y)|,
\end{align*}
where $\{\zeta_j;  j_1 \leq 0 < j_2\} \subseteq \mS$ is a bounded set. For all $m_1, m_2 \in \N$, we thus have
\begin{align*}
    |E_{(0,j_2)} (x,y)|\lesssim \sum_{j_1 \leq 0} \prod_{\mu =1 ,2} 2^{j_{\mu}Q_{\mu}}  (1+2^{j_{\mu}}|y_{\mu}^{-1}x_{\mu}|_{\mu})^{-m_{\mu}}.
\end{align*}
As before, $(1+2^{j_1}|y_1^{-1}x_1|_1)^{-m_1} \lesssim (1+|y_1^{-1}x_1|_1)^{-m_1}$ for $j_1 \leq 0$ on a bounded set $\Omega_1$. The Schwartz kernel $E_{(0,j_2)}(x, y)$ thus satisfies the desired estimate:
\beqq
    |E_{(0,j_2)} (x,y)| \lesssim   (1+|y_1^{-1}x_1|_1)^{-m_1} 2^{j_2Q_2} (1+2^{j_2}|y_2^{-1}x_2|_2)^{-m_2}.
\eeqq
Before proving that $\big\{(E_{(0, j_2)}, 2^{-(0,j_2)});  j_2 >0\big\}$ satisfies the second condition \eqref{BSEOs X}, we observe that by Remark \ref{remark pull out derivatives} and the decomposition \eqref{decomp PK}, we can ``pull out derivatives'' and write
\begin{equation}\label{eq Dj notation}
\begin{split}
    D_j = \Op(\zeta_j^{(2^j)}) &= \sum_{|\ap_2|= |\beta_2| =1} (2^{-j_2} X^2)^{\ap_2} \Op(\zeta_{j,\ap_2,\beta_2}^{(2^j)})(2^{-j_2} X^2)^{\beta_2},\\
    & =: \sum_{|\ap_2|= |\beta_2| =1} (2^{-j_2} X^2)^{\ap_2}D_{j, \ap_2, \beta_2} (2^{-j_2} X^2)^{\beta_2},
\end{split}
\end{equation}
where $\big\{\zeta_{j,\ap_2,\beta_2}; |\ap_2|= |\beta_2|=1\big\} \subseteq \mS_0^{\{1,2\}}$ is bounded. It suffices to prove \eqref{BSEOs X} with differential operators $(2^{-j}  X^2)^{\ap}$ on the left. The proof of \eqref{BSEOs X} with differential operators $(2^{-j} X^2)^{\beta}$ on the right is similar. Using the notation above in \eqref{eq Dj notation}, we write
\begin{align*}
    E_{(0,j_2)} = \sum_{j_1 \leq 0} \eta [D_j,\psi_1] \eta' = \sum_{j_1 \leq 0} \sum_{|\ap_2|=1} \eta [(2^{-j_2} X^2)^{\ap_2}D_{j,\ap_2},\psi_1] \eta'.
\end{align*}
We commute the differential operators $X^2\in \mathfrak{g}_2$ with the test function $\psi_1 \in C^{\infty}_0(\R^{q_1})$. As such, 
\begin{align*}
    E_{(0, j_2)}= \sum_{j_1 \leq 0} \sum_{|\ap_2|=1}\eta \lp (2^{-j_2}X^2)^{\ap_2}[D_{j, \ap_2}, \psi_1] \rp \eta'.
\end{align*}
By the product rule, 
\begin{align*}
    E_{(0, j_2)} = \sum_{j_1 \leq 0} \sum_{|\ap_2|=1} \eta  (2^{-j_2}X^2)^{\ap_2}[D_{j,\ap_2}, \psi_1] \eta' + 2^{-j_2 \deg \ap_2} \eta [D_{j,\ap_2}, \psi_1]  \widetilde{\eta}',
\end{align*}
where $(2^{-j_2}X^2)^{\ap_2} \eta' = 2^{-j_2\deg \ap_2}\widetilde{\eta}' \in C^{\infty}_0(\Omega)$. By the product rule again,
\begin{align*}
    E_{(0,j_2)} =& \sum_{|\ap_2|=1} \sum_{j_1 \leq 0}    (2^{-j_2}X^2)^{\ap_2}\eta[D_{j,\ap_2}, \psi_1] \eta'  + 2^{-j_2\deg \ap_2} \wte [D_{j,\ap_2}, \psi_1] \eta' \\
    &+ 2^{-j_2\deg \ap_2} \eta [D_{j,\ap_2}, \psi_1]  \widetilde{\eta}',
\end{align*}
where $[\eta, (2^{-j_2}X^2)^{\ap_2}] = 2^{-j_2\deg \ap_2}\wte $, for some $\wte  \in C^{\infty}_0(\Omega)$. We have thus shown that $E_{(0,j_2)}$ is a sum of derivatives of operators of the same form as $E_{(0,j_2)}$. The set
\begin{align*}
    \Big\{ \Big(\sum_{j_1 \leq 0}\eta[D_{j,\ap_2}, \psi_1] \eta', 2^{-(0,j_2)} \Big) ,  \Big( 2^{-j_2(\deg \ap_2-1)}\sum_{j_1 \leq 0} \wte [D_{j,\ap_2}, \psi_1] \eta', 2^{-(0,j_2)} \Big) ,\\  \Big( 2^{-j_2(\deg \ap_2-1)} \sum_{j_1 \leq 0} \eta [D_{j,\ap_2}, \psi_1]  \widetilde{\eta}', 2^{-(0,j_2)} \Big) ;  j_2 >0\Big\}
\end{align*}
is thus a bounded set of elementary operators\footnote{If $\{(F_j, 2^{-j});  j>0\}$ is a bounded set of elementary operators, then so is $\{(2^{-jn}F_j, 2^{-j});  j>0\}$ for $n\geq 0$.}. 

We proceed to show that the third term in \eqref{four sums PK} corresponds to a scaled sum of $2^{-(j_1, 0)}$ elementary operators which we denote $E_{(j_1, 0)}$; that is, let
\beqq
     \sum_{0< j_1} 2^{-j_1} E_{(j_1,0)}:= \sum_{0<j_1} 2^{-j_1} \Bigg( 2^{j_1} \sum_{j_2 \leq 0}\eta [D_j, \psi_1] \eta'  \Bigg).
\eeqq
We first verify that $\{( E_{(j_1,0)}, 2^{-(j_1,0)}); j_1>0  \}$ satisfies condition \eqref{BSPEOs X}. By the mean value theorem, and by the boundedness of the set $\{\zeta_j;  j_2 \leq 0 < j_1\} \subseteq \mS$, for all $m_1, m_2 \in \N$,
\begin{align*}
    |E_{(j_1,0)} (x,y)| \lesssim 2^{j_1} \sum_{j_2 \leq 0}   |y_1^{-1}x_1|_1 \prod_{\mu =1,2}  2^{j_{\mu} Q_{\mu}}(1+2^{j_{\mu}}|y_{\mu}^{-1}x_{\mu}|_{\mu})^{-m_{\mu}}.
\end{align*}
By the boundedness of $\Omega_2$, for $j_2 \leq 0$, we have $(1+2^{j_2}|y_2^{-1}x_2|_2)^{-m_2} \lesssim (1+|y_2^{-1}x_2|_2)^{-m_2}$. The Schwartz kernel $E_{(j_1,0)}(x, y)$ thus satisfies:
\beqq
    |E_{(j_1,0)} (x,y)| \lesssim  2^{j_1Q_1} (1+2^{j_1}|y_1^{-1}x_1|_1)^{-m_1+1}(1+|y_2^{-1}x_2|_2)^{-m_2}.
\eeqq
By taking $m_1$ large enough, we obtain the desired estimate. 

To verify that the set $\big\{(E_{(j_1,0)}, 2^{-(j_1,0)});  j_1>0\big\}$ satisfies \eqref{BSEOs X}, observe that by Remark \ref{remark pull out derivatives}, for every $j \in \Z^2$,
\begin{equation}\label{eq Dj notation 2}
    D_j = \Op(\zeta_j^{(2^j)}) = \sum_{|\ap_1|= |\beta_1| =1} (2^{-j_1} X^1)^{\ap_1} \Op(\zeta_{j,\ap_1,\beta_1}^{(2^j)})(2^{-j_1} X^1)^{\beta_1},
\end{equation}
where $\big\{\zeta_{j,\ap_1,\beta_1};  |\ap_1|= |\beta_1|=1 \big\} \subseteq \mS_0^{\{1,2\}}$ is bounded. We again define $D_{j, \ap_1, \beta_1} := \Op(\zeta_{j,\ap_1,\beta_1}^{(2^j)})$. It suffices to detail the proof of \eqref{BSEOs X} with differential operators on the left. By ``pulling out derivatives'' on the left, we have
\begin{align*}
    E_{(j_1,0)} &= 2^{j_1}\sum_{j_2 \leq 0}  \eta [D_j, \psi_1] \eta'\\
    &= 2^{j_1} \sum_{j_2 \leq 0} \sum_{|\ap_1|=1} \eta \lp \lp (2^{-j_1}X^1)^{\ap_1}D_{j,\ap_1} \rp \psi_1 - \psi_1 \lp (2^{-j_1}X^1)^{\ap_1}D_{j,\ap_1} \rp  \rp \eta'.
\end{align*}
\vspace{-.1in}
By the product rule,
\begin{align*}
    E_{(j_1,0)}=2^{j_1} \sum_{j_2 \leq 0} \sum_{|\ap_1|=1} \eta \lp (2^{-j_1}X^1)^{\ap_1}[D_{j,\ap_1},\psi_1] \rp \eta' - 2^{-j_1\deg \ap_1 }\eta [ D_{j,\ap_1}, \wtp_1]  \eta',
\end{align*}
where $(2^{-j_1}X^1)^{\ap_1}\psi_1 = 2^{-j_1\deg \ap_1}\wtp_1$, for some $\wtp_1 \in C^{\infty}_0( \R^{q_1} )$. By the product rule and letting $(2^{-j_1}X^1)^{\ap_1} \eta' = 2^{-j_1\deg \ap_1}\wte' \in C^{\infty}_0(\Omega)$, we have
\begin{align*}
    E_{(j_1,0)} =& 2^{j_1} \sum_{j_2 \leq 0} \sum_{|\ap_1|=1} \eta \lp (2^{-j_1}X^1)^{\ap_1}[D_{j,\ap_1},\psi_1] \eta' \rp - 2^{-j_1\deg \ap_1} \eta [D_{j,\ap_1},\psi_1]  \wte' \\
    &- 2^{-j_1\deg \ap_1}\eta [ D_{j,\ap_1}, \wtp_1]  \eta',
\end{align*}
Finally, commuting $(2^{-j_1}X^1)^{\ap_1}$ with $\eta$, the previous equation is
\begin{multline*}
    =  \sum_{|\ap_1|=1} \sum_{j_2 \leq 0} 2^{j_1} \lp (2^{-j_1}X^1)^{\ap_1}\eta[D_{j,\ap_1},\psi_1] \eta' \rp -  2^{-j_1\deg \ap_1}2^{j_1} \wte [D_{j,\ap_1},\psi_1] \eta'  \\
    - 2^{-j_1\deg \ap_1 } 2^{j_1} \eta [D_{j,\ap_1},\psi_1] \wte'  -2^{-j_1\deg \ap_1 }2^{j_1} \eta [ D_{j,\ap_1}, \wtp_1]  \eta',
\end{multline*}
where $[\eta, (2^{-j_1}X^1)^{\ap_1}] = 2^{-j_1\deg \ap_1} \wte \in C^{\infty}_0(\Omega)$. The set
\begin{align*}
    &\Big\{ \Big( 2^{j_1} \sum_{j_2 \leq 0} \eta  [D_{j,\ap_1},\psi_1] \eta' , 2^{-(j_1,0)} \Big), \Big( 2^{-j_1(\deg \ap_1 -1)} 2^{j_1}\sum_{j_2 \leq 0}\wte [D_{j,\ap_1},\psi_1] \eta',2^{-(j_1,0)} \Big), \\ 
    &\Big( 2^{-j_1(\deg \ap_1 -1)} 2^{j_1} \sum_{j_2 \leq 0}\eta [D_{j,\ap_1},\psi_1] \wte', 2^{-(j_1,0)} \Big),\hspace{-.05in}
    \Big( 2^{-j_1(\deg \ap_1 -1)} 2^{j_1}\sum_{j_2 \leq 0} \eta [ D_{j,\ap_1}, \wtp_1]  \eta',2^{-(j_1,0)} \Big); j_1>0 \Big\}
\end{align*}
is therefore a bounded set of elementary operators. 

It remains to show that the fourth and last term in \eqref{four sums PK} is a scaled sum of $2^{-(j_1, j_2)}$ elementary operators which we denote by $E_j$; that is, let
\beq\label{eq iv PK}
    \sum_{j_1, j_2>0} 2^{-j_1} E_j:= \sum_{\substack{j_1, j_2 > 0}} 2^{-j_1} \lp 2^{j_1} \eta [D_j, \psi_1] \eta' \rp .
\eeq
We first show that $\big\{\lp E_j, 2^{-j}\rp; j_1, j_2>0 \big\}$ satisfies \eqref{BSPEOs X}\footnote{As before, \eqref{BSPEOs X} for general $|\ap|, |\beta| \neq 0$ follows directly by Leibniz rule. So it suffices to consider the case $|\ap| = |\beta| = 0$. }. By the mean value theorem and by the boundedness of the set $\{\zeta_j;  j_1, j_2>0 \} \subseteq \mS$, 
\begin{align*}
    |E_j(x,y)| \lesssim 2^{j_{1} Q_1}(1+2^{j_1}|y_1^{-1}x_1|_1)^{-m_1+1} 2^{j_{2} Q_2}(1+2^{j_2}|y_2^{-1}x_2|_2)^{-m_2}, 
\end{align*}
where $m_1, m_2 \in \N$. Before verifying that $\big\{\lp E_j, 2^{-j}\rp; \ j_1, j_2>0 \big\}$ satisfies \eqref{BSEOs X}, note that by Remark \ref{remark pull out derivatives}, for every $j \in \Z^2$, we can ``pull out derivatives'' in both $\R^{q_1}$ and $\R^{q_2}$:
\begin{equation}\label{eq Dj notation 3}
\begin{split}
    D_j   &= \sum_{|\ap|= |\beta| =(1, 1)} (2^{-j} X)^{\ap} \Op(\zeta_{j,\ap,\beta}^{(2^j)})(2^{-j} X)^{\beta},
\end{split}
\end{equation}
where $\big\{\zeta_{j,\ap,\beta};  |\ap|, |\beta| =(1,1)\big\} \subseteq \mS_0^{\{1,2\}}$ is bounded. As before, we denote $D_{j, \ap, \beta}:=\Op(\zeta_{j,\ap,\beta}^{(2^j)})$. We can thus write
\begin{align*}
    E_j = 2^{j_1} \eta [D_j, \psi_1] \eta' = 2^{j_1}\sum_{|\ap|=(1,1)}  \eta [(2^{-j}X)^{\ap} D_{j,\ap} , \psi_1] \eta'.
\end{align*}
Notice that for $|\ap|=(1,1)$, $[(2^{-j}X)^{\ap}, \psi_1 ] \equiv 0$. We can thus pull the differential operator out of the commutator and write 
\begin{align*}
    E_j= 2^{j_1} \sum_{|\ap|=(1,1)}  \eta  \lp(2^{-j}X)^{\ap} [D_{j,\ap},  \psi_1] \rp  \eta'.
\end{align*}
By the product rule again,
\begin{align*}
    E_j = 2^{j_1} \sum_{|\ap|=(1,1)}  \eta (2^{-j}X)^{\ap} [D_{j,\ap},  \psi_1] \eta'   + 2^{-(j \cdot \deg \ap)} \eta [D_{j,\ap},  \psi_1] \wte',
\end{align*}
where $(2^{-j}X)^{\ap} \eta' = 2^{-j \cdot \deg \ap} \ \wte' \in C^{\infty}_0(\Omega)$. By commuting $\eta$ and $(2^{-j}X)^{\ap}$ and by applying the product rule, the previous equation is 
\begin{equation*}
    =2^{j_1}  \sum_{|\ap|=(1,1)} \hspace{-0.1in} (2^{-j}X)^{\ap}  \eta[D_{j,\ap},  \psi_1] \eta'   + 2^{-(j \cdot \deg \ap)} \Big( \wte [D_{j,\ap},  \psi_1] \eta'  
    +  \eta [D_{j,\ap},  \psi_1] \wte' \Big),
\end{equation*}
where $[\eta, (2^{-j}X)^{\ap}]= 2^{-(j \cdot \deg \ap)} \wte \in C^{\infty}_0(\Omega)$. Therefore, the set
\begin{align*}
    \Big\{\lp 2^{j_1}\eta[D_{j,\ap},  \psi_1] \eta', 2^{-j}\rp,  \lp 2^{-j \cdot \deg \ap + j \cdot (1, 1)}  2^{j_1}\wte [D_{j,\ap},  \psi_1] \eta', 2^{-j} \rp, \\
    \lp 2^{-j \cdot \deg \ap + j \cdot (1, 1)} 2^{j_1} \eta [D_{j,\ap},  \psi_1] \wte', 2^{-j}\rp;  j_1, j_2>0\Big\}
\end{align*}
is a bounded set of elementary operators\footnote{Indeed, $2^{-j\deg \ap + j\cdot (1,1)} \leq 1$, for all $j_1, j_2 >0$. }. Thus concluding the proof of Lemma \ref{lemma T psi1 pk}.
\end{proof}

\begin{proof}[Proof of Lemma \ref{comm SP}]
By commutativity of the differential operators $\mJ^1_{s_1}$, $\mJ^2_{s_2}$, Lemma 5.2 in \cite{CG84} implies estimate \eqref{comm J SP}. It thus remains to prove estimate \eqref{comm T SP}. By the growth condition \eqref{GC PK}, the Schwartz kernel $K(y^{-1} x)$ of $T$ is smooth away from the ``cross'': $ x_1 = y_1$ or $x_2 = y_2$. Let $\ep_1:= \frac{1}{4n_1!}$. By further localizing with $\phi_1 \prec \wt{\phi}_1$, where $\wt{\phi}_1 \in C^{\infty}_0(\R^{q_1})$, we write 
\begin{equation}\label{eq phi J1-ep and L psi}
    \begin{split}
        \norm{\phi_1 \mJ_{(\ep_1, 0)} [T, \psi_1]  f}_{L^2} \leq &  \norm{\phi_1 \mJ_{(\ep_1 -1, 0)} \wt{\phi}_1 (I+\mL_1) [T, \psi_1]  f}_{L^2} \\
        &+ \norm{\phi_1 \mJ_{(\ep_1 -1, 0)} (1- \wt{\phi}_1) (I+\mL_1) [T, \psi_1]  f}_{L^2}.
    \end{split}
\end{equation}
To bound the second term on the right-hand side of the inequality in \eqref{eq phi J1-ep and L psi}, note that the Schwartz kernel of the operator $\phi_1 \mJ_{1-\ep_1}^1 (1-\wt{\phi}_1)$ can be identified with a Schwartz function. By Young's inequality, 
\begin{align*}
    \norm{\phi_1 \mJ_{(\ep_1 -1, 0)} (1-\wt{\phi}_1) (I+\mL_1) [T, \psi_1]  f}_{L^2}  \lesssim &  \norm{J^1_{\ep_1-1}}_{L^1(G_1 \backslash \{0\})}  \norm{T}_{\mB(L^2)} \norm{ f}_{L^2}\\
    &+\norm{\phi_1 \mJ_{(\ep_1 -1, 0)} (1-\wt{\phi}_1) \mL_1 [T, \psi_1]  f}_{L^2}. 
\end{align*}
We can bound the second term on the right-hand side of the inequality above by Young's inequality. We thus obtain
\begin{equation*}
    \begin{split}
        \norm{\phi_1 \mJ_{(\ep_1 -1, 0)} (1-\wt{\phi}_1) (I+\mL_1) [T, \psi_1]  f}_{L^2} \lesssim & \norm{f}_{L^2}\\
        &+ \norm{\wt{\mL}_1((1- \vp_1) J^1_{\ep_1-1})}_{L^1(\R^{q_1})}\norm{T}_{\mB(L^2)} \norm{f}_{L^2}. 
    \end{split}
\end{equation*}
where $\vp_1 \in C^{\infty}_0(\R^{q_1})$ is supported near the identity so that $(1- \vp_1) J^1_{\ep_1-1} \in \mS$, and $\wt{\mL}_1$ is a left-invariant differential operator. 
To bound the first term on the right-hand side of the inequality in \eqref{eq phi J1-ep and L psi}, let $\wt{\phi}_1 \prec \doubletilde{\phi}_1 \in C^{\infty}_0(\R^{q_1})$. Since $I+\mL_1$ is a local operator, it remains to prove the $L^2$-boundedness of $\phi_1 \mJ_{(\ep_1 -1, 0)} \wt{\phi}_1 (I+\mL_1) \doubletilde{\phi}_1 [T, \psi_1] $ on $C^{\infty}_0(\R^q)$. 

It suffices to consider the localized operator 
\begin{align*}
    \phi_1 \mJ_{(\ep_1 -1, 0)} \wt{\phi}_1 (I+\mL_1) \doubletilde{\phi}_1 [T, \psi_1] \g_1,
\end{align*}
where $\doubletilde{\phi}_1 \prec \g_1 \in C^{\infty}_0(\R^{q_1})$. Indeed, the operator $\phi_1 \mJ_{(\ep_1 -1, 0)} \wt{\phi}_1 (I+\mL_1) \doubletilde{\phi}_1 [T, \psi_1] (1-\g_1)$ is $L^2$-bounded as $\supp \doubletilde{\phi}_1 \cap \supp (1- \g_1) = \emptyset$. This follows by a straightforward application of the Cotlar-Stein lemma. 


Observe that $\mJ_{\ep_1 -1}^1$ is given by group convolution on $\R^{q_1}$ with a Calderón-Zygmund kernel of order $1-4n_1!$ (see Proposition 5.1 in \cite{CG84}). In addition, $(I+\mL_1)$ is a Calderón-Zygmund operator of order $4n_1!$ acting on $C^{\infty}_0(\R^{q_1})$. By the algebra property stated in Theorem \ref{PSIO algebra}, the composition $\phi_1  \mJ_{(\ep_1 -1, 0)} \wt{\phi}_1 \circ (I+\mL_1) \doubletilde{\phi}_1$ is a Calderón-Zygmund operator of order $1$ on $C^{\infty}_0(\Omega_1)$, for some open, relatively compact set $\Omega_1 \Subset \R^{q_1}$. 

Consider a partition of unity on $\R^{q_2}$ of the form $1 = \sum_j \phi_2^j$, where $\phi_2^j \in C^{\infty}_0(\Omega_2^j)$, for a countable family of open, relatively compact sets $\Omega_2^j \Subset \R^{q_2}$ with the finite intersection property. 

By Lemma \ref{lemma T psi1 pk}, for $\doubletilde{\phi}_1 \prec \tripletilde{\phi}_1 \prec \g_1$, the operator $\tripletilde{\phi}_1 \otimes \phi_2^j[T,\psi_1] \g_1 \otimes \phi_2^k: C^{\infty}_0(\Omega_1 \times \Omega_2^{j,k}) \rightarrow C^{\infty}( \Omega_1 \times \Omega_2^{j,k})$ is a product singular integral operator of order $(-1, 0)$, where $\Omega_2^j, \Omega_2^k \subset \Omega_2^{j,k} \Subset \R^{q_2}$ and $\Omega_2^{j,k}$ is an open, relatively compact set. 

By the algebra property in Theorem \ref{PSIO algebra}, the operator 
\begin{align*}
    \phi_1  \mJ_{(\ep_1 -1, 0)} \wt{\phi}_1  \circ  (I+\mL_1)  \doubletilde{\phi}_1 \circ \tripletilde{\phi}_1 \otimes \phi_2^j [T, \psi_1] \g_1 \otimes \phi_2^k
\end{align*}
is a product singular integral operator mapping $C^{\infty}_0(\Omega_1 \times \Omega_2^{j,k}) \rightarrow C^{\infty}( \Omega_1 \times \Omega_2^{j,k})$ of order $(0,0)$. It remains to show that the sum 
\begin{align}\label{eq sum on whole space}
    \sum_{j,k} \phi_1 \mJ_{(\ep_1 -1, 0)} \wt{\phi}_1 (I+\mL_1) \doubletilde{\phi}_1 \otimes \phi_2^j [T, \psi_1] \g_1 \otimes \phi_2^k 
\end{align} 
is a product singular integral operator $C^{\infty}_0(\R^q) \rightarrow C^{\infty} (\R^q)$ of order $(0,0)$ on the whole space. 

The growth condition as stated in Definition \ref{def PSIO} follows directly from the growth condition for each summand and the finite intersection property of the sets $\Omega_2^j$.

It thus remains to verify that \eqref{eq sum on whole space} satisfies the cancellation condition as stated in Definition \ref{def PSIO}. For all bounded sets of bump functions $\mC_{\mu} = \{(\zeta_{\mu}, z_{\mu}, r_{\mu})\} \subseteq C^{\infty}_0(\R^{q_{\mu}}) \times \R^{q_{\mu}} \times (0, \infty)$, where $\mu= 1, 2$, we define the operator $\Big( \sum_{j,k} \phi_1 \mJ_{(\ep_1 -1, 0)} \wt{\phi}_1 (I+\mL_1) \doubletilde{\phi}_1 \otimes \phi_2^j [T, \psi_1] \g_1 \otimes \phi_2^k \Big)^{\zeta_1, x_1}$ by
\begin{align*}
    &\int_{\R^{q_1}} \Big\la \Big( \sum_{j,k} \phi_1 \mJ_{(\ep_1 -1, 0)} \wt{\phi}_1 (I+\mL_1) \doubletilde{\phi}_1 \otimes \phi_2^j [T, \psi_1] \g_1 \otimes \phi_2^k \Big)^{\zeta_1, x_1} (\zeta_2), \vp_2 \Big\ra \vp_1(x_1) dx_1\\
    &= \Big\la \sum_{j,k} \phi_1 \mJ_{(\ep_1 -1, 0)} \wt{\phi}_1 (I+\mL_1) \doubletilde{\phi}_1 \otimes \phi_2^j [T, \psi_1] \g_1 \otimes \phi_2^k(\zeta_1 \otimes \zeta_2), \vp_1 \otimes \vp_2 \Big\ra. 
\end{align*}

By the definition of product singular integral operators, each summand composed with dilated vector fields
\begin{align}\label{eq CZO zeta1 x1}
    r_1^{Q_1} (r_1 X_{x_1}^1)^{\ap} \lp \phi_1 \mJ_{(\ep_1 -1, 0)} \wt{\phi}_1 (I+\mL_1) \doubletilde{\phi}_1 \otimes \phi_2^j [T, \psi_1] \g_1 \otimes \phi_2^k \rp^{\zeta_1, x_1}
\end{align}
is a Calderón-Zygmund operator $C^{\infty}_0( \Omega_2^{j,k}) \rightarrow C^{\infty}(\Omega_2^{j,k})$. To be more precise, for each multi-index $\ap$, the Calderón-Zygmund operator in \eqref{eq CZO zeta1 x1} satisfies the growth condition and the cancellation condition uniformly over the bounded set of bump functions $\mC_1 \subseteq C^{\infty}_0(\R^{q_1}) \times \R^{q_1} \times (0, \infty) $ given above and $x_1 \in \R^{q_1}$. As such, the sum of operators
\begin{align}\label{eq sum op}
    r_1^{Q_1} (r_1 X_{x_1}^1)^{\ap} \lp \sum_{j,k} \phi_1 \mJ_{(\ep_1 -1, 0)} \wt{\phi}_1 (I+\mL_1) \doubletilde{\phi}_1 \otimes \phi_2^j [T, \psi_1] \g_1 \otimes \phi_2^k \rp^{\zeta_1, x_1}
\end{align}
satisfies the growth condition for Calderón-Zygmund operators on $\R^{q_2}$ by the finite intersection property of the sets $\Omega_2^j$. 


We need to verify that the sum \eqref{eq sum op} satisfies the cancellation condition for Calderón-Zygmund operators on $\R^{q_2}$. Given a bounded set of bump functions $\mC_2 =\{ (\zeta_2, z_2, r_2)\} \subseteq C^{\infty}_0(\R^{q_2}) \times \R^{q_2} \times (0, \infty)$, for every multi-index $\beta$, we need to bound
\begin{align*}
    \sup_{\substack{(\zeta_2, z_2, r_2) \in \mC_2;  \\ x_2 \in \R^{q_2}}}  \Bigg|r_2^{Q_2}(r_2X_{x_2}^2)^{\beta} \Bigg( r_1^{Q_1} (r_1 X_{x_1}^1)^{\ap} \Bigg( \sum_{j,k} \phi_1 \mJ_{(\ep_1 -1, 0)} \wt{\phi}_1 (I+\mL_1) \doubletilde{\phi}_1 \otimes \phi_2^j [T, \psi_1] \g_1 \otimes \phi_2^k \Bigg)^{\zeta_1, x_1} \Bigg) \\ (\zeta_2)  (x_2) \Bigg|. 
\end{align*}
Observe that the expression above is
\begin{align*}
    =\sup_{\substack{(\zeta_2, z_2, r_2) \in \mC_2; \\ x_2 \in \R^{q_2}}}  \Bigg| r_1^{Q_1} (r_1 X_{x_1}^1)^{\ap} \Bigg(r_2^{Q_2} (r_2X^2_{x_2})^{\beta}  \Bigg( \phi_1 \mJ_{(\ep_1 -1, 0)} \wt{\phi}_1 (I+\mL_1) \doubletilde{\phi}_1  [T, \psi_1] \g_1   \Bigg)^{\zeta_2, x_2} \Bigg) (\zeta_1)  (x_1) \Bigg|.
\end{align*}
In fact, the right-hand side of the equality above is
\begin{align}\label{eq Tchi2 x2}
    =\sup_{\substack{(\zeta_2, z_2, r_2) \in \mC_2;  \\ x_2 \in \R^{q_2}}}  \Bigg| r_1^{Q_1} (r_1 X_{x_1}^1)^{\ap}  \Bigg( \phi_1 \mJ_{(\ep_1 -1, 0)} \wt{\phi}_1 (I+\mL_1) \doubletilde{\phi}_1  \Big[r_2^{Q_2} (r_2X^2_{x_2})^{\beta}T^{\zeta_2, x_2}, \psi_1 \Big] \g_1   \Bigg) (\zeta_1)  (x_1) \Bigg|. 
\end{align}
Recall that $T: C^{\infty}_0(\R^q) \rightarrow C^{\infty}(\R^q)$ is a product singular integral operator of order $(0, 0)$. As such, the operator
\begin{align*}
    r_2^{Q_2} (r_2 X^2_{x_2})^{\beta}T^{\zeta_2, x_2}
\end{align*}
is a Calderón-Zygmund operator of order $0$ on $\R^{q_1}$ with bounds uniform in $\mC_2$, and $x_2$, for every multi-index $\beta$. 

On the one hand, by the single-parameter Calderón-Zygmund theory, we know that the operator $  [r_2^{Q_2} (r_2X_{x_2}^2)^{\beta}T^{\zeta_2, x_2}, \psi_1] $ is a Calderón-Zygmund operator of order $-1$. By further localizing with test functions $\doubletilde{\phi}_1, \g_1 \in C^{\infty}_0(\R^{q_1})$, we obtain a Calderón-Zygmund operator of order $-1$ on a relatively compact open set. 

On the other hand, as noted above, by the algebra property stated in Theorem \ref{PSIO algebra}, the composition $\phi_1  \mJ_{(\ep_1 -1, 0)} \wt{\phi}_1 \circ (I+\mL_1) \doubletilde{\phi}_1$ is a Calderón-Zygmund operator of order $1$ on an open, relatively compact set.

By the algebra property in Theorem \ref{PSIO algebra}, we can thus conclude that the operator
\begin{align}
     \phi_1 \mJ_{(\ep_1 -1, 0)} \wt{\phi}_1 (I+\mL_1) \doubletilde{\phi}_1  [r_2^{Q_2} (r_2 X^2_{x_2})^{\beta}T^{\zeta_2, x_2}, \psi_1] \g_1
\end{align}
is a Calderón-Zygmund operator of order $0$ on an open, relatively compact set with bounds uniform in $\mC_2$ and $x_2$, for every multi-index $\beta$.

By the cancellation condition for Calderón-Zygmund operators of order $0$, we conclude that \eqref{eq Tchi2 x2} is uniformly bounded, as desired. 

Following the argument above, we can show that the operator 
\begin{align}\label{eq sum op 2}
    r_2^{Q_2} (r_2 X_{x_2}^2)^{\ap}\Big( \sum_{j,k} \phi_1 \mJ_{(\ep_1 -1, 0)} \wt{\phi}_1 (I+\mL_1) \doubletilde{\phi}_1 \otimes \phi_2^j [T, \psi_1] \g_1 \otimes \phi_2^k \Big)^{\zeta_2, x_2}
\end{align}
is a Calderón-Zygmund operator on $\R^{q_1}$ satisfying the growth condition and the cancellation condition uniformly over the bounded set of bump functions  $\mC_2 \subseteq C^{\infty}_0(\R^{q_2}) \times \R^{q_2} \times (0, \infty) $ given above and $x_2 \in \R^{q_2}$.

By retracing the proof, we can show that the transpose of the operator \eqref{eq sum on whole space} satisfies the cancellation condition for product singular integral operators of order $(0, 0)$ on $\R^{q_1} \times \R^{q_2} $.

We can thus conclude that the operator 
\begin{align*}
    \sum_{j,k} \phi_1 \mJ_{(\ep_1 -1, 0)} \wt{\phi}_1 (I+\mL_1) \doubletilde{\phi}_1 \otimes \phi_2^j [T, \psi_1] \g_1 \otimes \phi_2^k
\end{align*}
is a product singular integral operator mapping $C^{\infty}_0(\R^q) \rightarrow C^{\infty}(\R^q)$ of order $(0,0)$. By Theorem \ref{PSIO L2 bdd}, it follows that the operator is $L^2$-bounded as desired. 
\end{proof}

We record an intermediate single-parameter a priori estimate next.

\begin{lemma}\label{CG lemma SP}
If $T$ is $L^2$-invertible, then for all $l_1 \in \N$, $\psi_1 \in C^{\infty}_0(\R^{q_1})$, and $f \in C^{\infty}_0(\R^q)$, 
\begin{equation}\label{key estimate SP}
    \RLs{\psi_1 f}{l_1 \ep_1,0} \leq C(l_1, \psi_1, T)  \big(\RLs{\psi_1 Tf}{l_1\ep_1,0}+ \RLs{ f}{l_1 \ep_1 - \ep_1,0} \big), 
\end{equation}
where the implicit constant depends on $T$ and $T^{-1}$ in an admissible way\footnote{The implicit constant depends on $T$ and $T^{-1}$ in an admissible way as defined in Proposition \ref{CG lemma PK}. }.
\end{lemma}

\begin{proof}
    For $\psi_1 \prec \eta_1 \in C^{\infty}_0(\R^{q_1})$, one can show that the operator $(1- \eta_1) \mL_1  T \psi_1$ is $L^2$-bounded. This follows from a straightforward application of the Cotlar-Stein lemma and the observation that $\supp (1- \eta_1) \cap \supp \psi_1 = \emptyset$. The proof of of Lemma \ref{CG lemma SP} then follows by applying this observation along with the multi-parameter commutator estimates in Lemma \ref{comm SP} to the proof of Lemma 5.3 in \cite{CG84}, to which we refer the reader. 
\end{proof}

\subsection{Proof of the A Priori Estimate for Product Kernels} \ 

We record the base case of Proposition \ref{CG lemma PK} in the following lemma.

\begin{lemma}\label{lemma base case PK} For all $\psi_{\mu} \prec \eta_{\mu} \in C^{\infty}_0(\R^{q_{\mu}})$, where $\mu =1, 2$, 
\begin{equation*}
    \RLs{\psi_1 \otimes \psi_2 f}{\ep_1,\ep_2} \lesssim    \RLs{\psi_1 \otimes \psi_2 Tf}{\ep_1,\ep_2}\hspace{-0.06in}+ \RLs{\eta_1 T f}{\ep_1,0}\hspace{-0.06in}  + \RLs{\eta_2 T f}{0,\ep_2} \hspace{-0.07in} +\|f\|_{L^2}, 
\end{equation*}
where $f\in C^{\infty}_0(\R^q)$ and where the implicit constant depends on the listed cutoff functions and the operators $T, T^{-1}$ in an admissible way\footnote{See the definition of an \textit{admissible constant} in Proposition \ref{CG lemma PK}. }. 
\end{lemma}

\begin{proof}[Proof of Lemma \ref{lemma base case PK}]
    Let $\ep_1:= \frac{1}{4n_1!},   \ep_2:= \frac{1}{4n_2!}$. By applying Lemma \ref{CG lemma SP} to $\mJ_{(0, \ep_2)}\psi_2 f$, we have
\begin{align*}
    \RLs{\psi_1 \mJ_{(0, \ep_2)}\psi_2 f}{0, \ep_1} \lesssim \RLs{\psi_1 T \mJ_{(0, \ep_2)}\psi_2 f}{\ep_1,0} + \|\mJ_{(0, \ep_2)}\psi_2 f\|_{L^2}.
\end{align*}
By commuting left-invariant and right-invariant operators and introducing commutators for operators that do not commute, we obtain
\begin{equation}\label{eq 1 lemma base case PK}
\begin{split}
    \RLs{\psi_1 \otimes \psi_2 f}{\ep_1,\ep_2}  \lesssim & \RLs{\psi_1  \otimes \psi_2 T  f}{\ep_1, \ep_2}  + \RLs{\psi_1 \mJ_{(0,\ep_2)}  [T,\psi_2] f}{\ep_1, 0} \\
    &+ \RLs{\psi_2 f}{0, \ep_2}.
    \end{split}
\end{equation}
It remains to bound the last two terms on the right-hand side of the inequality above. For the second to last term in \eqref{eq 1 lemma base case PK}, let $\psi_2 \prec \phi_2 \in C^{\infty}_0(\R^{q_2})$ and $\psi_1 \prec \eta_1$. By the triangle inequality,
\begin{equation}\label{eq in proof lemma 3.21}
    \begin{split}
    \RLs{\psi_1 \mJ_{(0,\ep_2)}  [T,\psi_2] f}{\ep_1, 0} \leq & \RLs{\psi_1 \otimes \phi_2 \mJ_{(0,\ep_2)}  [T,\psi_2] \eta_1 f}{\ep_1, 0} \\
    &+ \RLs{\psi_1 \otimes \phi_2 \mJ_{(0,\ep_2)}  [T,\psi_2] (1-\eta_1) f}{\ep_1, 0} \\
    &+\RLs{\psi_1 \otimes(1-\phi_2) \mJ_{(0,\ep_2)}  [T,\psi_2] \eta_1 f}{\ep_1, 0} \\
    &+\RLs{\psi_1 \otimes(1-\phi_2) \mJ_{(0,\ep_2)}  [T,\psi_2] (1-\eta_1 ) f}{\ep_1, 0}.
    \end{split}
\end{equation}
By applying Lemma \ref{comm SP} to the first term on the right-hand side of the inequality in \eqref{eq in proof lemma 3.21}, 
\begin{align*}
     \RLs{\psi_1 \otimes \phi_2 \mJ_{(0,\ep_2)}  [T,\psi_2] \eta_1 f}{\ep_1, 0} \lesssim \RLs{\eta_1 f}{\ep_1, 0}. 
\end{align*}
We bound the second term on the right-hand side of the inequality in \eqref{eq in proof lemma 3.21} by applying Lemma \ref{comm SP} and noting that $J_{\ep_1-1}^1 \in L^1(G_1)$.
\begin{equation}
    \begin{split}
        &\RLs{\psi_1 \otimes \phi_2 \mJ_{(0,\ep_2)}  [T,\psi_2] (1-\eta_1) f}{\ep_1, 0} \\
        &\hspace{1in}\lesssim \norm{J_{\ep_1-1}^1}_{L^1(G_1)} (\norm{f}_{L^2} + \norm{  \phi_2 \mJ_{(0,\ep_2)}  [\mL_1\psi_1 T(1-\eta_1),\psi_2]  f}_{L^2}). 
    \end{split}
\end{equation}
$T$ is localized away from its singularity in $\R^{q_1}$. As such, $\phi_2 \mJ_{(0,\ep_2)}  [\mL_1\psi_1 T(1-\eta_1),\psi_2]$ is an $L^2$-bounded operator. This can be shown by retracing the proof of \eqref{comm T SP} after replacing $T$ with the $L^2$-bounded operator $\mL_1\psi_1 T(1-\eta_1)$. We thus have
\begin{equation*}
    \RLs{\psi_1 \otimes \phi_2 \mJ_{(0,\ep_2)}  [T,\psi_2] (1-\eta_1) f}{\ep_1, 0} \lesssim \norm{f}_{L^2}. 
\end{equation*}
To bound the third term on the right-hand side of the inequality in \eqref{eq in proof lemma 3.21}, let $\psi_2 \prec \g_2 \prec \phi_2$, 
\begin{align*}
    & \hspace{-.5in} \RLs{\psi_1 \otimes(1-\phi_2) \mJ_{(0,\ep_2)}  [T,\psi_2] \eta_1 f}{\ep_1, 0} \\
     \leq & \RLs{\psi_1 \otimes(1-\phi_2) \mJ_{(0,\ep_2 -1)} \g_2 (I+\mL_2)  [T,\psi_2] \eta_1 f}{\ep_1, 0} \\
    & + \RLs{\psi_1 \otimes(1-\phi_2) \mJ_{(0,\ep_2 -1)} (1-\g_2) (I+\mL_2)  [T,\psi_2] \eta_1 f}{\ep_1, 0}. 
\end{align*}
By associativity of convolution, since $J_{\ep_2-1 }^2 \in \mS(G_2 \backslash \{0\})$, the right-hand side of the inequality above is
\begin{align*}
    \lesssim &  \norm{ ((1- \chi_2 ) J_{\ep_2 -1}^2) *_2 (I+\mL_2)\de_0}_{L^1(G_2)}\RLs{ \psi_1 [T,\psi_2] \eta_1 f}{\ep_1, 0} \\
    &+ \norm{J_{\ep_2 -1}^2}_{L^1(G_2)} \RLs{ (1-\g_2) (I+\mL_2)  [T,\psi_2] \eta_1 f}{\ep_1, 0}, 
\end{align*}
where $\supp \chi_2 \subset \{|x_2|_2 \leq 1\}$. The right-hand side of the inequality above in turn is
\begin{align*}
    & \lesssim \norm{T}_{\mB(L^2)} \RLs{ \eta_1 f}{\ep_1, 0} + \RLs{  (I+\mL_2)(1-\wt{\g}_2)  T \psi_2 \eta_1 f}{\ep_1, 0},
\end{align*}
for some $\psi_2 \prec \wt{\g}_2 \in C^{\infty}_0(\R^{q_2})$. By Cotlar-Stein's lemma, one can show that the operator $(I+\mL_2) (1- \wt{\g}_2) T \psi_2$ is $L^2$-bounded. The third term in \eqref{eq in proof lemma 3.21} is thus
\begin{align*}
    \RLs{\psi_1 \otimes(1-\phi_2) \mJ_{(0,\ep_2)}  [T,\psi_2] \eta_1 f}{\ep_1, 0} & \lesssim \RLs{ \eta_1 f}{\ep_1, 0} + \norm{f}_{L^2}.
\end{align*}
To bound the fourth and last term in \eqref{eq in proof lemma 3.21}, let $\psi_2 \prec \g_2 \prec \phi_2$. By the triangle inequality,
\begin{align*}
    & \hspace{-1in} \RLs{\psi_1 \otimes(1-\phi_2) \mJ_{(0,\ep_2)}  [T,\psi_2] (1-\eta_1 ) f}{\ep_1, 0} \\
    \leq & \RLs{\psi_1 \otimes(1-\phi_2) \mJ_{(0,\ep_2-1)} \g_2 (I+\mL_2)  [T,\psi_2] (1-\eta_1 ) f}{\ep_1, 0}\\
    &+\RLs{(1-\phi_2) \mJ_{(0,\ep_2-1)} \psi_1 \otimes(1-\g_2) (I+\mL_2)  [T,\psi_2] (1-\eta_1 ) f}{\ep_1, 0}. 
\end{align*}
Recalling that, away from the identity in $\R^{q_{\mu}}$, $J_{\ep_{\mu}-1}^{\mu}\in L^1(G_{\mu})$ are Schwartz, the right-hand side of the inequality above is
\begin{align*}
    \leq & \norm{J_{\ep_1-1}}_{L^1(G_1)} \norm{ (1-\chi_2) J_{(0,\ep_2-1)} *_2 (I+\mL_2) \de_0}_{L^1(G_2)} \norm{(I+\mL_1)\psi_1  [T,\psi_2] (1-\eta_1 ) f}_{ L^2 }\\
    &+ \norm{J_{\ep_1-1}}_{L^1(G_1)} \norm{J_{\ep_2 -1}^2}_{L^1(G_2)} \norm{ (I+\mL_1) (I+\mL_2) \psi_1 \otimes (1-\wt{\g}_2)  T  (1-\eta_1 ) \otimes \psi_2 f}_{L^2}. 
\end{align*}
By Cotlar-Stein's lemma, one can show that $(I+\mL_1) \psi_1 [T, \psi_2] (1- \eta_1)$, and $(I+\mL_1) (I+\mL_2) \psi_1 \otimes (1-\wt{\g}_2)  T  (1-\eta_1 ) \otimes \psi_2$ are both $L^2$-bounded.
The second to last term in \eqref{eq 1 lemma base case PK} is thus
\begin{equation}
    \begin{split}
        \RLs{\psi_1 \mJ_{(0,\ep_2)}  [T,\psi_2] f}{\ep_1, 0} & \lesssim \RLs{\eta_1 f}{\ep_1, 0} + \norm{f}_{L^2}.
    \end{split}
\end{equation}
It thus remains to bound $\RLs{\psi_2 f}{0, \ep_2}$ and $\RLs{\eta_1 f}{\ep_1, 0}$. By symmetry, it suffices to estimate $\RLs{\psi_2 f}{0, \ep_2}$. By the $L^2$-invertibility of $T$, we have
\begin{align*}
    \RLs{\psi_2 f}{0, \ep_2} & \lesssim \RLs{ \psi_2 T f }{0, \ep_2} + \norm{\mJ_{(0, \ep_2)} [T, \psi_2] f }_{L^2}. 
\end{align*}
We proceed to bound the second term on the right-hand side of the inequality above. Let $\psi_2 \prec \phi_2$. By the triangle inequality, we have
\begin{align*}
    \norm{\mJ_{(0, \ep_2)} [T, \psi_2] f }_{L^2}  \leq & \norm{\phi_2 \mJ_{(0, \ep_2)} [T, \psi_2] f }_{L^2} + \norm{ (1-\phi_2) \mJ_{(0, \ep_2)} [T, \psi_2] f }_{L^2}. 
\end{align*}
By Lemma \ref{comm SP}, we bound the localized commutator term; that is,
\begin{align*}
    \norm{\mJ_{(0, \ep_2)} [T, \psi_2] f }_{L^2}  \lesssim & \norm{f}_{L^2} + \norm{ (1-\phi_2) \mJ_{(0, \ep_2)} [T, \psi_2] f }_{L^2}. 
\end{align*}
It thus remains to bound $\norm{ (1-\phi_2) \mJ_{(0, \ep_2)} [T, \psi_2] f }_{L^2}$. Let $\psi_2 \prec \g_2 \prec \phi_2$. By the triangle inequality,
\begin{align*}
    \norm{ (1-\phi_2) \mJ_{(0, \ep_2)} [T, \psi_2] f }_{L^2} \leq & \norm{ (1-\phi_2) \mJ_{(0, \ep_2 -1)} \g_2 (I+\mL_2) [T, \psi_2] f }_{L^2} \\
    &+ \norm{ (1-\phi_2) \mJ_{(0, \ep_2 -1)} (1-\g_2) (I+\mL_2) [T, \psi_2] f }_{L^2}. 
\end{align*}
As above, we obtain
\begin{align*}
    \norm{ (1-\phi_2) \mJ_{(0, \ep_2)} [T, \psi_2] f }_{L^2} \leq & \norm{(1-\chi_2) J_{\ep_2 -1}^2 *_2 (I+\mL_2)\de_0 }_{L^1} \norm{  [T, \psi_2] f }_{L^2} \\
    &+ \norm{J_{\ep_2-1}^2}_{L^1(G_2)} \norm{ (I+\mL_2)  (1-\wt{\g}_2)  T \psi_2 f }_{L^2}, 
\end{align*}
for some test function $\psi_2 \prec \wt{\g}_2 \in C^{\infty}_0(\R^{q_2})$. By the Cotlar-Stein lemma, the operator $ (I+\mL_2)  (1-\wt{\g}_2)  T \psi_2$ is $L^2$-bounded. Thus concluding the proof of Lemma \ref{lemma base case PK}. 
\end{proof}

\begin{proof}[Proof of Proposition \ref{CG lemma PK}] By commuting the test functions and the differential operators on each factor space, 
\begin{equation*}
    \begin{split}
        \RLs{  \psi_1  \otimes \psi_2 f}{l_1\ep_1, l_2\ep_2} \leq & \RLs{\psi_1 \otimes \psi_2 \mJ_{(l_1 \ep_1-\ep_1, l_2 \ep_2-\ep_2)} f}{\ep_1,\ep_2}\\
        &+ \RLs{ [\psi_1, \mJ_{(l_1 \ep_1-\ep_1, 0)}]\otimes [\psi_2, \mJ_{(0, l_2\ep_2-\ep_2)}]f}{\ep_1,\ep_2}.
    \end{split}
\end{equation*}
By \eqref{comm J SP} applied to the second term on the right-hand side of the inequality above, we have 
\begin{equation*}
    \begin{split}
        \RLs{ \psi_1  \otimes \psi_2 f}{l_1\ep_1, l_2\ep_2} \lesssim & \RLs{\psi_1 \otimes \psi_2 \mJ_{(l_1 \ep_1-\ep_1, l_2 \ep_2-\ep_2)} f}{\ep_1,\ep_2}\\
        &+ \RLs{ f}{l_1 \ep_1-\ep_1, l_2\ep_2-\ep_2}. 
    \end{split}
\end{equation*}
Let $\psi_{\mu} \prec \eta_{\mu}$, for $\mu =1, 2$. It remains to bound the first term on the right-hand side of the inequality above. By the base case in Lemma \ref{lemma base case PK},
\begin{align*}
    \RLs{\psi_1 \otimes \psi_2 \mJ_{(l_1 \ep_1 -\ep_1, l_2 \ep_2 -\ep_2)} f}{\ep_1,\ep_2}  \lesssim & \RLs{\psi_1 \otimes \psi_2 T\mJ_{(l_1 \ep_1-\ep_1, l_2 \ep_2-\ep_2)} f}{\ep_1,\ep_2} \\
    & + \RLs{\eta_1 T \mJ_{(l_1 \ep_1-\ep_1, l_2 \ep_2-\ep_2)} f}{\ep_1,0} \\
    &+ \RLs{\eta_2 T \mJ_{(l_1 \ep_1-\ep_1, l_2 \ep_2-\ep_2)} f}{0,\ep_2} + \|\mJ_{(l_1 \ep_1-\ep_1, l_2 \ep_2-\ep_2)}f\|_{L^2}. 
\end{align*}
By commuting left- and right-invariant operators, the right-hand side of the above is
\begin{align*}
    = & \RLs{\psi_1 \otimes \psi_2 \mJ_{(l_1 \ep_1-\ep_1, l_2 \ep_2-\ep_2)} T f}{\ep_1,\ep_2}  + \RLs{\eta_1 \mJ_{(l_1 \ep_1-\ep_1, 0)} T f}{\ep_1,l_2 \ep_2-\ep_2}\\
    &+ \RLs{\eta_2  \mJ_{(0, l_2 \ep_2-\ep_2)} T f}{l_1 \ep_1-\ep_1, \ep_2}  +\RLs{f}{l_1\ep_1-\ep_1, l_2 \ep_2-\ep_2}.
\end{align*}
Introducing commutators for the cutoff functions and the operators $\mJ_{l_{\mu}\ep_{\mu}-\ep_{\mu}}^{\mu}$ for $\mu =1, 2$, the right-hand side of the equality above is
\begin{align*}
    \leq  \RLs{ \psi_1 \otimes \psi_2 T f}{l_1 \ep_1, l_2 \ep_2)} + \RLs{ [\psi_1, \mJ_{(l_1 \ep_1-\ep_1, 0)}] \otimes [\psi_2, \mJ_{(0, l_2 \ep_2-\ep_2)}] T f}{\ep_1,\ep_2} \\
    +  \RLs{\eta_1 T f}{l_1\ep_1,l_2 \ep_2-\ep_2} + \RLs{[\eta_1, \mJ_{(l_1 \ep_1-\ep_1, 0)}] T f}{\ep_1,l_2 \ep_2-\ep_2}\\
    + \RLs{\eta_2  T f}{l_1 \ep_1-\ep_1, l_2 \ep_2}+ \RLs{[\eta_2 , \mJ_{(0, l_2 \ep_2-\ep_2)}] T f}{l_1 \ep_1-\ep_1, \ep_2} \\
    +\RLs{f}{l_1\ep_1-\ep_1, l_2 \ep_2-\ep_2}.
\end{align*}
By the commutator estimates in Lemma \ref{comm SP}, the equation above is
\begin{align*}
    \lesssim&  \RLs{ \psi_1 \otimes \psi_2 T f}{l_1 \ep_1, l_2 \ep_2)} 
    +  \RLs{\eta_1 T f}{l_1\ep_1,l_2 \ep_2-\ep_2} \\
    &+ \RLs{\eta_2  T f}{l_1 \ep_1-\ep_1, l_2 \ep_2}+ \RLs{T f}{l_1 \ep_1-\ep_1, l_2 \ep_2-\ep_2} +\RLs{f}{l_1\ep_1-\ep_1, l_2 \ep_2-\ep_2}.
\end{align*}
After commuting the right-invariant operators and $T$, the desired estimate follows. Thus concluding the proof of Proposition \ref{CG lemma PK}. 
\end{proof}

\subsection{Proof of the Inversion Theorem for Product Kernels} \ 

The non-isotropic Sobolev norms and the usual Euclidean Sobolev norms are related by the following estimate which follows from a straightforward adaptation of Proposition 5.1.27 p.274 in \cite{Str14} to which we refer the reader for the proof. 
\begin{proposition}\label{sob MP} For $k \in \N$, and $f \in C^{\infty}_0(\Omega)$, where $\Omega$ is a relatively compact open subset of $\R^{q_1} \times \R^{q_2}$, 
\begin{equation}
    \|f\|_{L^2_k} \lesssim  \RLs{f}{\frac{k}{4(n_1-1)!}, \frac{k}{4(n_2-1)!}}.
\end{equation}
\end{proposition}

Suppose $T$ is invertible on $L^2$ with bounded inverse $T^{-1}$. By the Schwartz kernel theorem and the left-translation invariance of $T$, we know that $T^{-1}$ is also given by $T^{-1}g = g*L$, for some distribution $L \in \mD'(\R^q)$. To prove Theorem \ref{main thm PK} for product kernels, we need to verify that $L$ is in fact a product kernel. To do so, we first establish the following lemma. 

\begin{lemma}\label{lemma L smooth}
$L(t_1, t_2) \in C^{\infty}(\R^{q_1} \times \R^{q_2} \backslash \{t_1 =0\} \cup \{t_2 =0\} )$. 
\end{lemma}
\begin{proof}[Proof of Lemma \ref{lemma L smooth}]
    Choose $\ap_1, \ap_2>0$ so that $J_{(-\ap_1, -\ap_2)} \in L^2(G)$ and $J^{\mu}_{-\ap_{\mu}} \in \mS(\R^{q_{\mu}} \backslash 0)$ (see Proposition 5.1 in \cite{CG84}). By the $L^2$-boundedness of $T^{-1}$, we have that $J_{(-\ap_1, -\ap_2)} * L \in L^2$. It remains to show that $J_{(-\ap_1, -\ap_2)} * L \in C^{\infty}(\R^{q_1} \times \R^{q_2} \backslash \{t_1 =0\} \cup \{t_2 =0\})$. 

    It will then follow that $J_{(\ap_1, \ap_2)} * J_{(-\ap_1, -\ap_2)} * L = L$, in the sense of distributions, is also in $C^{\infty}(\R^{q_1} \times \R^{q_2} \backslash \{t_1 =0\} \cup \{t_2 =0\})$ (see Proposition 5.1 in \cite{CG84}). 

    Let $\Omega_1 \times \Omega_2 \Subset \R^{q_1} \times \R^{q_2}$ be an open relatively compact set s.t. $0 \notin \overline{\Omega}_1$ and $0 \notin \overline{\Omega}_2$. Let $\{\chi_{\eta}\}_{\eta>0}$ be an approximation of the identity on $\R^q$ and $\phi_{\mu}^j \in C^{\infty}_0(\R^{q_{\mu}})$, $\phi_{\mu}^j \equiv 1$ on $\overline{\Omega}_{\mu}$,
    $\phi_{\mu}^j \equiv 0$ near $0$,  for $\mu =1, 2$ and $\phi_{\mu}^j \prec \eta_{\mu}^j \prec \phi_{\mu}^{j+1}$ for $j \in \N$. 
For $s_1, s_2>0$, by Proposition \ref{CG lemma PK}, 
\begin{equation}\label{eq proof lemma 3.24}
\begin{split}
    &\hspace{-1in}\RLs{\phi_1^1 \otimes \phi_2^1 (\phi_1^2 \otimes \phi_2^2 (J_{(-\ap_1, -\ap_2)} * L) * \chi_{\eta})}{s_1, s_2} \\
    \lesssim & \RLs{\phi_1^1 \otimes \phi_2^1 T (\phi_1^2 \otimes \phi_2^2 (J_{(-\ap_1, -\ap_2)} * L) * \chi_{\eta})}{s_1, s_2}  \\
    &+ \RLs{\eta^1_1  T (\phi_1^2 \otimes \phi_2^2 (J_{(-\ap_1, -\ap_2)} * L) * \chi_{\eta})}{s_1, s_2-\ep_2} \\
    &+ \RLs{\eta_2^1  T (\phi_1^2 \otimes \phi_2^2 (J_{(-\ap_1, -\ap_2)} * L) * \chi_{\eta})}{s_1 -\ep_1, s_2} \\
    &+ \RLs{\phi_1^2 \otimes \phi_2^2 (J_{(-\ap_1, -\ap_2)} * L) * \chi_{\eta}}{s_1 -\ep_1, s_2 -\ep_2}.
    \end{split}
\end{equation}
We need to show that all four terms on the right-hand side of the inequality \eqref{eq proof lemma 3.24} are finite.

To bound the first term on the right-hand side of the inequality \eqref{eq proof lemma 3.24}, first note that $T(J_{(-\ap_1, -\ap_2)} * L) = J_{(-\ap_1, - \ap_2)}$ in the sense of distributions. As such, $\phi_1^1 \otimes \phi_2^1 T(J_{(-\ap_1, -\ap_2)} * L) \in \mS$. We introduce more cutoff functions and write:
\begin{equation}\label{eq 1st term}
\begin{split}
    &\hspace{-1in} \RLs{\phi_1^1 \otimes \phi_2^1 T (\phi_1^2 \otimes \phi_2^2 (J_{(-\ap_1, -\ap_2)} * L) * \chi_{\eta})}{s_1, s_2} \\
    \leq & \RLs{\phi_1^1 \otimes \phi_2^1 T(J_{(-\ap_1, -\ap_2)} * L * \chi_{\eta}) }{s_1, s_2} \\
    &+ \RLs{\phi_1^1 \otimes \phi_2^1 T(\phi_1^2 \otimes (1-\phi_2^2)(J_{(-\ap_1, -\ap_2)} * L* \chi_{\eta}))}{s_1, s_2}\\
    &+ \RLs{\phi_1^1 \otimes \phi_2^1 T((1-\phi_1^2) \otimes \phi_2^2(J_{(-\ap_1, -\ap_2)} * L* \chi_{\eta}))}{s_1, s_2}\\
    &+ \RLs{\phi_1^1 \otimes \phi_2^1 T((1-\phi_1^2) \otimes (1-\phi_2^2)(J_{(-\ap_1, -\ap_2)} * L* \chi_{\eta}))}{s_1, s_2}.
    \end{split}
\end{equation}
As $\eta \rightarrow 0$, the first term on the right-hand side of the inequality converges to \\
$\RLs{\phi_1^1 \otimes \phi_2^1 J_{(-\ap_1, -\ap_2)}}{s_1, s_2} < \infty$. 
It thus remains to bound the latter three terms on the right-hand side of the inequality \eqref{eq 1st term}. By symmetry, we bound the second and third terms similarly. We will thus only detail the proof for the second term. The operator $\mJ_{(0, s_2)}\phi_2^1 T (1- \phi_2^2)$ is $L^2$-bounded\footnote{This follows from the Cotlar-Stein lemma, and a straightforward adaptation of Lemma 1.1.19 in \cite{Str14} to the setting of graded Lie groups paired with the observation that $\supp \phi_2^1 \cap \supp (1-\phi_2^2)=\emptyset$.}. We thus have
\begin{align*}
    &\RLs{\phi_1^1 \otimes \phi_2^1 T(\phi_1^2 \otimes (1-\phi_2^2)(J_{(-\ap_1, -\ap_2)} * L* \chi_{\eta}))}{s_1, s_2}\lesssim \RLs{\phi_1^2 (J_{(-\ap_1, -\ap_2)} * L* \chi_{\eta})}{s_1, 0}. 
\end{align*}
Let $\phi_1^2 \prec \phi_1^3$. By Lemma \ref{CG lemma SP}, we have
\begin{align*}
    \RLs{   \phi_1^2 \phi_1^3 (J_{(-\ap_1, -\ap_2)} * L* \chi_{\eta})}{s_1, 0} \lesssim & \RLs{\phi_1^2 T \phi_1^3 (J_{(-\ap_1, -\ap_2)} * L* \chi_{\eta})}{s_1, 0} \\
    &+ \RLs{\phi_1^3 (J_{(-\ap_1, -\ap_2)} * L* \chi_{\eta})}{s_1-\ep_1, 0}.
\end{align*}
By the triangle inequality, we bound the right-hand side of the inequality above by
\begin{align*}
    \leq & \RLs{\phi_1^2 T  (J_{(-\ap_1, -\ap_2)} * L* \chi_{\eta})}{s_1, 0} + \RLs{\phi_1^2 T (1-\phi_1^3) (J_{(-\ap_1, -\ap_2)} * L* \chi_{\eta})}{s_1, 0} \\
    &+ \RLs{\phi_1^3 (J_{(-\ap_1, -\ap_2)} * L* \chi_{\eta})}{s_1-\ep_1, 0}.
\end{align*}
Observe that the operator $\mJ_{(s_1, 0)}\phi_1^2 T(1-\phi_1^3)$ is $L^2$ bounded. It thus remains to bound the last term on the right-hand side of the inequality above. Repeating this process with $\phi_1^3 \prec \phi_1^4 \prec \cdots \prec \phi_1^N$, we obtain that the second term on the right-hand side of the inequality in \eqref{eq 1st term} is bounded. 

Finally, for the fourth and final term on the right-hand side of the inequality in \eqref{eq 1st term}, we have 
\begin{align*}
    &\hspace{-.5in} \RLs{\phi_1^1 \otimes \phi_2^1 T((1-\phi_1^2) \otimes (1-\phi_2^2)(J_{(-\ap_1, -\ap_2)} * L* \chi_{\eta}))}{s_1, s_2}\\
    & \lesssim \sum_{|(\beta_1, \beta_2)| \leq (k,k)} \norm{ p_{(\beta_1, \beta_2)} (X^{(\beta_1, \beta_2)}K)}_{L^1(G \backslash \{t_1=0\} \cup \{t_2 =0\})}  \norm{J_{(-\ap_1, -\ap_2)} * L * \chi_{\eta}}_{L^2}, 
\end{align*}
for some $k >0$, homogeneous polynomials $p_{(\beta_1, \beta_2)}$ and left-invariant differential operators $X^{(\beta_1, \beta_2)}$. Indeed, in this case, the product kernel $K$ can be identified with a decaying smooth function. 

By symmetry, we bound the second and third terms on the right-hand side of the inequality \eqref{eq proof lemma 3.24} similarly. We thus only detail the proof for the second summand. We further localize with $\phi_1^1 \prec \eta_1^1 \prec \phi_1^2$. By the triangle inequality,
\begin{align*}
    &\hspace{-.5in} \RLs{\eta^1_1  T (\phi_1^2 \otimes \phi_2^2 (J_{(-\ap_1, -\ap_2)} * L) * \chi_{\eta})}{s_1, s_2-\ep_2} \\
    \leq & \RLs{\eta^1_1  T ( \phi_2^2 (J_{(-\ap_1, -\ap_2)} * L) * \chi_{\eta})}{s_1, s_2-\ep_2} \\
    &+ \RLs{\eta^1_1  T ((1-\phi_1^2) \otimes \phi_2^2 (J_{(-\ap_1, -\ap_2)} * L) * \chi_{\eta})}{s_1, s_2-\ep_2}.
\end{align*}
The operator $\mJ_{(s_1, 0)} \eta^1_1  T (1-\phi_1^2)$ is $L^2$-bounded as noted above. The right-hand side of the inequality above is thus
\begin{equation}\label{eq 2nd term before commutator}
\begin{split}
    \lesssim & \RLs{\eta^1_1  T ( \phi_2^2 (J_{(-\ap_1, -\ap_2)} * L) * \chi_{\eta})}{s_1, s_2-\ep_2} + \RLs{  \phi_2^2 (J_{(-\ap_1, -\ap_2)} * L) * \chi_{\eta}}{0, s_2-\ep_2}.
    \end{split}
\end{equation}
By repeatedly applying the single-parameter a priori estimate with a sequence of cutoff functions $\phi_2^2 \prec \phi_2^3 \prec \cdots \prec \phi_2^N$, we bound the second term in \eqref{eq 2nd term before commutator}. Let $\phi_2^2 \prec \phi_2^3$ and $\eta_1^1 \prec \eta_1^2$. By localizing further and by the triangle inequality, we bound the first term on the right-hand side of the inequality in \eqref{eq 2nd term before commutator} as follows. 
\begin{align*}
    &\hspace{-.5in} \RLs{\eta^1_1  T ( \phi_2^2 (J_{(-\ap_1, -\ap_2)} * L) * \chi_{\eta})}{s_1, s_2-\ep_2} \\
    \leq & \RLs{\eta^1_1 \otimes \phi_2^3 T ( \eta_1^2 \otimes \phi_2^2 (J_{(-\ap_1, -\ap_2)} * L) * \chi_{\eta})}{s_1, s_2-\ep_2}\\
    &+\RLs{\eta^1_1 \otimes \phi_2^3 T ( (1-\eta_1^2) \otimes \phi_2^2 (J_{(-\ap_1, -\ap_2)} * L) * \chi_{\eta})}{s_1, s_2-\ep_2}\\
    &+\RLs{\eta^1_1 \otimes (1-\phi_2^3) T (\eta_1^2 \otimes \phi_2^2 (J_{(-\ap_1, -\ap_2)} * L) * \chi_{\eta})}{s_1, s_2-\ep_2}\\
    &+\RLs{\eta^1_1 \otimes (1-\phi_2^3) T ((1-\eta_1^2)  \otimes \phi_2^2 (J_{(-\ap_1, -\ap_2)} * L) * \chi_{\eta})}{s_1, s_2-\ep_2}.
\end{align*}
Observe that all four terms can be bounded by using ideas detailed above. 

Finally, we reapply the a priori estimate to the last term in \eqref{eq proof lemma 3.24}. By repeatedly following this procedure using a sequence of cutoff functions $\phi_{\mu}^j \prec \phi_{\mu}^{j+1}$ and finally taking $\eta \rightarrow 0$, we obtain the desired result. 
\end{proof}

\begin{proof}[Proof of Theorem \ref{main thm PK} for product kernels]
We need to verify that $L$ satisfies both the growth and cancellation conditions for product kernels. By Lemma \ref{lemma L smooth} and by scaling considerations, the proof of \eqref{GC PK} reduces to proving that:
\begin{equation}\label{eq reduction GC PK}
\sup_{|t_1|_1, |t_2|_2 \sim 1} |\p_{t_1}^{\ap_1} \p_{t_2}^{\ap_2} L(t_1, t_2)| \lesssim 1.
\end{equation}
Indeed, given \eqref{eq reduction GC PK}, for all $R_1, R_2>0$, the kernel $R_1^{Q_1} R_2^{Q_2} L(R_1 \cdot  t_1, R_2 \cdot t_2)$, is the kernel associated to the operator $D_{(R_1, R_2)} T^{-1} D_{(R_1, R_2)}^{-1}$, where we define $D_{(R_1, R_2)} f(x_1, x_2) := f(R_1 \cdot x_1, R_2 \cdot x_2)$. The admissible constants for the operators $D_{(R_1, R_2)} T^{-1} D_{(R_1, R_2)}^{-1}$ in the a priori estimate are uniformly bounded in $R_1, R_2>0$. Hence, for all $x_1, x_2 \neq 0$, by writing $x_1= R_1 \cdot t_1$ and $x_2 = R_2 \cdot t_2$ for $|t_1|_1, |t_2|_2 \sim 1$, we obtain
\begin{align*}
    |\p_{x_1}^{\ap_1} \p_{x_2}^{\ap_2} L(x_1, x_2)| \lesssim |x_1|_1^{-Q_1 - \deg \ap_1} |x_2|_2^{-Q_2 - \deg \ap_2}. 
\end{align*}
It thus remains to show that the growth condition holds for $L$ restricted to $|t_1|_1, |t_2|_2 \sim 1$. Let $\phi_1 \otimes \phi_2 \in C^{\infty}_0(\R^{q_1} \backslash \{0\}) \otimes  C^{\infty}_0(\R^{q_2} \backslash \{0\} )$ such that $\supp \phi_1 \otimes \phi_2(t_1, t_2) \equiv 1$ on $\{ |t_1|_1, |t_2|_2 \sim 1 \}$ and such that $L$ and each of its derivatives, up to some finite order $m$ chosen below, do not change signs on $\supp \phi_1 \otimes \phi_2$. By the Sobolev embedding, 
\begin{equation*}
    \sup_{|t_1|_1,|t_2|_2 \sim 1}|\p_{t_1}^{\ap_1} \p_{t_2}^{\ap_2} L(t_1,t_2)| \lesssim \|\phi_1 \otimes \phi_2 L\|_{L^2_m},
\end{equation*}
for some $m>0$. By the choice of the cutoff functions, 
\begin{align}
    \norm{\phi_1 \otimes \phi_2 L}_{L^2_m} \leq \sum_{|\ap| \leq m} \lf \int  \p^{\ap_1}_{x_1} \p_{x_2}^{\ap_2} \Big[\phi_1(x_1) \phi_2(x_2) L(x_1, x_2) \Big] dx \rf.
\end{align}
On the support of $\phi_1 \otimes  \phi_2$, the distribution $L$ can be identified with a smooth function and each of its derivatives introduced above do not change signs. In this way, we can choose $\psi_{\mu}$ with support comparable to $\phi_{\mu}$ near the annulus $|x_{\mu}|_{\mu } \sim 1$ and choose $\zeta_{\mu}$ to have support in a very small open set about the origin so that the above equation is 
\begin{align}
    \leq  \sum_{|\ap| \leq m} \lf  \int \int \p^{\ap_1}_{x_1} \p^{\ap_2}_{x_2} \Big[ \psi(x_1) \psi(x_2) L(y_1^{-1}x_1, y_2^{-1}x_2) \Big]  \zeta_1(y_1) \zeta_2(y_2) dy dx \rf.
\end{align}
We rewrite the expression above as follows:
\begin{equation*}
    \sum_{|\ap| \leq m} \lf  \int \p^{\ap_1}_{x_1} \p^{\ap_2}_{x_2} \Big[  \psi_1(x_1) \psi_2(x_2) T^{-1} (\zeta_1\otimes \zeta_2)(x_1, x_2) \Big] dx \rf. 
\end{equation*}
By compactness, the expression above is 
\begin{align*}
    & \leq \norm{\psi_1 \otimes \psi_2 T^{-1} \zeta_1\otimes \zeta_2}_{L^2_m}.
\end{align*}
By Proposition \ref{sob MP}, there exist $s_1, s_2>0$,
\begin{equation*}
    \norm{\psi_1 \otimes \psi_2 T^{-1} \zeta_1\otimes \zeta_2}_{L^2_m} \leq \RLs{\psi_1 \otimes \psi_2 T^{-1} \zeta_1 \otimes \zeta_2}{s_1,s_2}.
\end{equation*}
The boundedness of the expression above follows from the right-invariance of the differential operators and the left-invariance of $T^{-1}$. Thus proving \eqref{eq reduction GC PK}.

In the next step of the proof of Theorem \ref{main thm PK}, we need to show that $L$ satisfies the cancellation condition \eqref{CC2 PK}. That is, we need to show that for $R_1>0$, a bounded set $\mB_1 \subseteq C^{\infty}_0(\R^{q_1})$, and $\phi_1 \in \mB_1$, the distribution $L_{\phi_1,R_1} \in C^{\infty}_0(\R^{q_2})'$, defined by
\begin{align*}
    L_{\phi_1,R_1} (t_2) = \int L(t_1, t_2) \phi_1(R_1 \cdot t_1) dt_1,
\end{align*}
is a Calderón-Zygmund kernel, with seminorms uniformly bounded in $\phi_1$ and $R_1$. $L_{\phi_2,R_2}(t_1)$ defined analogously must also correspond to a Calderón-Zygmund kernel. By symmetry, we only present the proof for $L_{\phi_1,R_1} (t_2)$. 

By homogeneity, we first prove that $L$ satisfies \eqref{CC2 PK} with $R_1 = 1$. By making use of the scale-invariant property of Calderón-Zygmund kernels\footnote{If $K(t_1)$ is a Calderón-Zygmund kernel on $\R^{q_1}$, then $R_1^{Q_1} K(R_1 t_1)$ is too. Moreover, their respective seminorms as defined in \eqref{GC PK} and \eqref{seminorm2 PK} are equal.}, proving that $L_{\phi_1,1}$ satisfies the growth condition \eqref{GC PK} reduces to proving the following estimate. For all $\ap_2 \in \N^{q_2}$, 
\begin{equation}\label{eq GC on annuli}
    \sup_{|t_2|_2 \sim 1} |\p_{t_2}^{\ap_2}  L_{\phi_1, 1}(t_2)| \lesssim 1.
\end{equation}
We pick $\phi_2 \in C^{\infty}_0(\R^{q_2})$ s.t. $\phi_1 \equiv 1$ for $|t_2|_2 \sim 1$ s.t. $L_{\phi_1, 1}(t_2)$ and all of its derivatives up to some finite order $m_2$ determined below, do not change signs on $\supp \phi_2$. By the Sobolev embedding, there exists $m_2 \in \N$ s.t.
\begin{align*}
    \sup_{|t_2|_2 \sim 1} |\p_{t_2}^{\ap_2}  L_{\phi_1, 1}(t_2)| \lesssim \norm{\phi_2 L_{\phi_1, 1} }_{L^2_{m_2}(\R^{q_2})}.
\end{align*}
By the choice of cutoff function, the right-hand side of the inequality above is 
\begin{align*}
    \leq \sum_{|\ap_2| \leq m_2} \lf \int \p^{\ap_2}_{x_2} \phi_1(x_1) \Big [ \phi_2(x_2) L(x_1, x_2) \Big] dx  \rf.
\end{align*}
There exists $\psi_2, \zeta_2 \in C^{\infty}_0(\R^{q_2})$ with $\supp \psi_2 \cap \supp \zeta_2 = \emptyset$ and $\psi_1 \in C^{\infty}_0(\R^{q_1})$ with $0 \in \supp \psi_1$ s.t. the right-hand side of the inequality above is 
\begin{align*}
    & \leq \sum_{|\ap_2| \leq m_2} \lf \int \p^{\ap_2}_{x_2} \psi_1(x_1)  \Big[ \psi_2(x_2) T^{-1} (\wt{\phi}_1 \otimes \zeta_2)(x_1, x_2)  \Big] dx  \rf,
\end{align*}
where $\wt{\phi}_1 (x_1) = \phi_1(x_1^{-1})$. By compactness followed by Proposition \ref{sob MP}, there exists $s_2>0$ s.t. the above expression is
\begin{align*}
    \lesssim \RLs{\psi_1 \otimes \psi_2 T^{-1} (\wt{\phi}_1 \otimes \psi_2)}{0, s_2}.
\end{align*}
Finally since the right-invariant differential operators commute with the left-invariant operator $T^{-1}$, we obtain the desired bound. The general growth condition for $L_{\phi_1,1}(t_2)$ follows directly by homogeneity considerations as described earlier.

We then need to show that $L_{\phi_1,1}$ satisfies the cancellation condition \eqref{CC1 PK} for Calderón-Zygmund kernels. That is, given a bounded set $\mB_2 \subseteq C^{\infty}_0(\R^{q_2})$ and $R_2>0$, we need to show that
\beqq
    \sup_{\phi_2 \in \mB_2; \ R_2>0} \lf \int L_{\phi_1,1} (t_2) \phi_2(R_2\cdot  t_2) dt_2 \rf \leq C_{\mB_2}.
\eeqq
By a standard scaling argument, it suffices to prove the cancellation condition holds for $R_2 = 1$. Let $\psi_1 \in C_0^{\infty}(\R^{q_1})$, $ \psi_2 \in C_0^{\infty}(\R^{q_2})$ be chosen s.t. $\psi_1 \otimes \psi_2(0, 0)=1$. By the Sobolev embedding followed by Proposition \ref{sob MP}, there exists some $ s_1, s_2>0$ s.t.
\begin{align*}
    \sup_{\phi_2 \in \mB_2} \lf \int L_{\phi_1, 1}(t_2) \phi_2(t_2) dt_2 \rf 
    &\lesssim \RLs{\psi_1 \otimes \psi_2 T^{-1} \phi_1 \otimes \phi_2}{s_1,s_2}.
\end{align*}
Finally, by commuting the right-invariant differential operators with the left-invariant operator $T^{-1}$, the right-hand side of the inequality above is bounded. 
Thus concluding the proof of the cancellation condition for $L_{\phi_1,1}$. We can thus in turn conclude that $L_{\phi_1,1}$ and $L_{\phi_2,1}$ are Calderón-Zygmund kernels on $\R^{q_2}$ and $\R^{q_1}$ respectively.

To conclude the proof of Theorem \ref{main thm PK} for product kernels, it remains to show that $L_{\phi_1,R_1}(t_2)$ is a Calder\'on-Zygmund kernel on $\R^{q_2}$ for $R_1>0$. To make use of the scale invariance of the operators at play, recall the following dilation operators $D_{(R_1, R_2)} f(x_1, x_2):= f(R_1\cdot x_1, R_2 \cdot x_2)$. 

For $R_1>0$, after a change of variables, we have
\begin{align*}
    |\p_{t_2}^{\ap_2} L_{\phi_1,R_1}(t_2)| &=   \Big|  \int R_1^{-Q_1} (\p_{t_2}^{\ap_2}L)(R_1^{-1}\cdot  t_1, t_2 ) \phi_1(t_1) dt_1 \Big|, 
\end{align*}
where $L^{R_1} (\cdot, \cdot) := R_1^{-Q_1} L(R_1^{-1}\cdot, \cdot)$ is the convolution kernel of $D_{(R_1, 1)}^{-1}  T^{-1} D_{(R_1, 1)}$. Observe that since $L_{\phi_1,1}$ is a Calderón-Zygmund kernel, the kernel $(L^{R_1})_{\phi_1,1}$ associated to the operator $D_{(R_1, 1)}^{-1}  T^{-1}   D_{(R_1, 1)}$ satisfies the following estimate:
\begin{align*}
    |\p_{t_2}^{\ap_2} (L^{R_1})_{\phi_1,1}(t_2)| & \lesssim |t_2|_2^{-Q_2 - \deg \ap_2},
\end{align*}
where the constant is uniform in $R_1$. $L_{\phi_1,R_1}$ thus satisfies the growth condition for Calderón-Zygmund kernels with seminorms uniformly bounded in $\phi_1 \in \mB_1$ and $R_1>0$:
\begin{align*}
    \sup_{\substack{\phi_1 \in \mB_1;\\
    R_1>0}}|\p_{t_2}^{\ap_2} L_{\phi_1,R_1}(t_2)| \lesssim |t_2|_2^{-Q_2 - \deg  \ap_2}. 
\end{align*}
Similarly, we can show that $L_{\phi_1,R_1}$ satisfies the cancellation condition for Calderón-Zygmund kernels with bounds independent of $\phi_1$ and $R_1$. To avoid redundancy, we omit this step and conclude the proof of Theorem \ref{main thm PK} in the case of product kernels. 
\end{proof}

\vspace{.2in}

\section{Inversion Theorem for Flag Kernels}

\begin{definition}\label{def FK}
A \textit{flag kernel} $K$ on $\R^q = \R^{q_1} \times \cdots \times \R^{q_{\nu}}$ is a distribution satisfying the following two conditions:

(i) Growth condition - For every multi-index $\ap = (\ap_1, \ldots, \ap_{\nu}) \in \N^{q_1} \times \cdots  \times \N^{q_{\nu}}$, there is a constant $C_{\ap}$ s.t.
\begin{equation}\label{GC FK}
|\p_{t_1}^{\ap_1} \cdots \p_{t_{\nu}}^{\ap_{\nu}} K(t)| \leq C_{\ap} \prod_{\mu=1}^{\nu} \lp |t_1|_1 + \ldots +|t_{\mu}|_{\mu} \rp^{-Q_{\mu} - \deg \ap_{\mu}}.
\end{equation}
We define the least possible $C_{\ap}$ to be a seminorm.

(ii) Cancellation condition - This condition is defined recursively. 

• For $\nu =1$, given a bounded set $\mB \subseteq C^{\infty}_0(\R^q)$, 
\begin{equation}\label{CC1 FK}
    \sup_{\phi \in \mB; \ R>0} \Big| \int K(t) \phi(R\cdot t) dt \Big| <\infty. 
\end{equation}

• For $\nu >1$, given $1 \leq \mu \leq \nu$, a bounded set $\mB_{\mu} \subseteq C^{\infty}_0(\R^{q_{\mu}})$, $\phi_{\mu} \in \mB_{\mu}$, and $R_{\mu} >0$, the distribution $K_{\phi_{\mu},R_{\mu}} $ defined by
\begin{equation}\label{CC2 FK}
    K_{\phi_{\mu}, R_{\mu}} (\ldots, t_{\mu -1},   t_{\mu+1}, \ldots) := \int K(t)\phi_{\mu}(R_{\mu} \cdot  t_{\mu})dt_{\mu}
\end{equation}
is a flag kernel on the $(\nu-1)$-factor space $\cdots \times \R^{q_{\mu-1}} \times \R^{q_{\mu+1}} \times \cdots $ where the bounds are independent of the choice of $\phi_{\mu}$ and $R_{\mu}$.

For the base case $\nu=0$, we define the space of flag kernels to be $\C$ with its usual topology. For every seminorm $|\cdot|$ on the space of $(\nu-1)$-factor flag kernels, we define a seminorm on flag kernels on $\R^{q_1} \times \cdots  \times \R^{q_{\nu}}$ by 
\beq\label{seminorm2 FK}
    |K|:=\sup_{\phi_{\mu} \in \mB_{\mu}, R_{\mu}>0}|K_{\phi_{\mu},R_{\mu}}|,
\eeq
which we assume to be finite.
\end{definition}
\begin{remark}
\cite{MRS95} and \cite{NRSW12} studied flag kernels on Heisenberg-type groups and on homogeneous groups respectively; while \cite{NRS01} studied flag kernels on $\nu$-factor product spaces and homogeneous groups. \cite{Glo10} and \cite{Glo13} investigated flag kernels on homogeneous groups independently. Other recent results on flag kernels include \cite{Yan09}, \cite{HLW19} and the references therein. 
\end{remark}

In an effort to highlight the main ideas of the proof, we again detail the $2$-parameter case. The general $\nu$-parameter case follows from a few straightforward modifications. Much like in the proof of Theorem \ref{main thm PK} for product kernels, the key idea in the proof of Theorem \ref{main thm PK} for flag kernels is an a priori estimate which we record in the next proposition. 

\begin{proposition}\label{CG lemma FK} There exist $\ep_1, \ep_2>0$ s.t. for all $l_{\mu}\in \N$, and $\psi_{\mu} \prec \eta_{\mu} \in C^{\infty}_0(\R^{q_{\mu}})$, where $\mu = 1,2$, for $f \in C^{\infty}_0(\R^q)$, 
\begin{equation*}
\begin{split}
    \RLs{\psi_1 \otimes \psi_2 f}{l_1 \ep_1, l_2 \ep_2} \lesssim & \RLs{\psi_1 \otimes \psi_2 Tf}{l_1 \ep_1, l_2 \ep_2} + \RLs{\eta_1 T f}{l_1\ep_1, l_2 \ep_2 - \ep_2} \\
    &+ \RLs{ f}{l_1 \ep_1 - \ep_1, l_2 \ep_2 -\ep_2},
\end{split}
\end{equation*}
where the implicit constant depends on $\psi_{\mu}, \eta_{\mu}$ and on the operators $T$ and $T^{-1}$ in an admissible\footnote{See the definition of an admissible constant in Proposition \ref{CG lemma PK}.} way. 
\end{proposition}
To avoid redundancy, we will only highlight the new ideas needed to adapt the proof of Proposition \ref{CG lemma PK} to that of Proposition \ref{CG lemma FK}. 


\begin{lemma}\label{lemma base case FK}
There exist $\ep_1, \ep_2>0$, s.t. for $\psi_{\mu} \in C^{\infty}_0(\R^{q_{\mu}})$, where $\mu =1,2$, and $\psi_1 \prec \eta_1$, for $f \in C^{\infty}_0(\R^q)$, 
\begin{equation*}
    \RLs{\psi_1 \otimes \psi_2 f}{\ep_1,\ep_2} \lesssim \RLs{\psi_1 \otimes \psi_2 Tf}{\ep_1,\ep_2} + \RLs{\eta_1 Tf}{\ep_1,0} + \norm{f}_{L^2},
\end{equation*}
where the implicit constant depends on $\psi_1, \psi_2, \eta_1$ and on the operators $T$ and $T^{-1}$ in an admissible\footnote{See the definition of an admissible constant in Proposition \ref{CG lemma PK}.} way. 
\end{lemma}
\begin{proof}[Proof of Lemma \ref{lemma base case FK}]
    As before, let $\ep_1 = \frac{1}{4n_1!}$ and $\ep_2 = \frac{1}{4n_2!}$. Applying the single-parameter a priori estimate in Lemma \ref{CG lemma SP} to $\mJ_{(\ep_1,0)} \psi_1 f$,
\begin{align*}
    \RLs{\psi_1 \otimes \psi_2 f}{\ep_1,\ep_2} \lesssim \norm{\mJ_{(\ep_1,\ep_2)}\psi_2 T\psi_1 f}_{L^2} +  \norm{\mJ_{(\ep_1,0)}\psi_1 f}_{L^2}.
\end{align*}
We introduce commutators for the first term and apply the single-parameter estimate in Lemma \ref{CG lemma SP} to the second term.
\begin{align*}
    \RLs{\psi_1 \otimes \psi_2 f}{\ep_1,\ep_2} &\lesssim \RLs{ \psi_1 \otimes \psi_2 Tf}{\ep_1, \ep_2} + \norm{\mJ_{(\ep_1,\ep_2)} \psi_2 [T,\psi_1] f}_{L^2}+ ( \RLs{\psi_1 Tf}{\ep_1,0} + \norm{f}_{L^2}).
\end{align*}
It remains to bound the second summand $\norm{\mJ_{(\ep_1,\ep_2)} \psi_2 [T,\psi_1] f}_{L^2}$. A flag kernel $K$ on a direct product of graded Lie groups is a product kernel presenting more singularity in the first variable. As such, we will show that not only is $[T,\psi_1]$ smoothing on $\R^{q_1}$ as shown in Lemma \ref{comm SP} but it is also smoothing on the second factor space $\R^{q_2}$. By Theorem \ref{PSIO theorem BSEOs}, the proof reduces to the proof of the following claim. 
\vspace{-.1in}
\begin{claim}\label{claim fk T psi1}
Given $\eta, \eta' \in C^{\infty}_0(\Omega)$, there exists a bounded set of elementary operators $\{(E_j, 2^{-j}); \ j \in \Z^2_{\geq 0} \}$ s.t.
\begin{equation}\label{eq [T,p] decomp flag}
    \eta [T,\psi_1] \eta' = \sum_{\substack{(j_1, j_2) \in \Z_{\geq 0}^2}} 2^{-j_2} E_{(j_1, j_2)}.
\end{equation}
\end{claim}
We have the following decomposition\footnote{See Lemma 4.2.24 in \cite{Str14} for the precise formulation or Corollary 2.4.4 in \cite{NRS01} for an analogous formulation.} of $T$:
\begin{equation*}
    T = \sum_{\substack{ j_1 \geq j_2}} \Op( \zeta_j^{(2^j)}) =:\sum_{\substack{ j_1 \geq j_2}} D_j,
\end{equation*}
where on the one hand $\{\zeta_j; \  j \in \Z^2, \ j_1 = j_2\}\subseteq \mS_0^{\{2\}}$ is bounded, and on the other hand, $\{\zeta_j; \  j \in \Z^2, \ j_1 >j_2\} \subseteq \mS_0^{\{1,2\}}$ is bounded. Consider the following four separate cases:
\begin{equation}\label{four sums FK}
\begin{split}
    \sum_{\substack{j;  j_1 \geq j_2}}  \eta  [D_j,\psi_1] \eta' 
    = &  \sum_{\substack{0 \geq j_1 \geq j_2}}  \ \eta [D_j,\psi_1] \eta'  + \sum_{\substack{ j_1 > 0 \geq j_2}}  \ \eta [D_j,\psi_1] \eta' \\
    &+ \sum_{\substack{ j_1 = j_2 > 0}} \eta  [D_j,\psi_1] \eta'+ \sum_{\substack{ j_1 > j_2 > 0}} \eta  [D_j,\psi_1] \eta'.
\end{split}
\end{equation}
The first term on the right-hand side above converges to an elementary operator $(E_0, 2^{0})$ as in the product kernel case (see \eqref{eq [T,p] decomp}). One can readily verify that the last term is a sum of elementary operators scaled by a factor $2^{-j_2}$ since $\{\zeta_j;  j_1>j_2>0 \}\subseteq \mS_0^{\{1,2\}}$ is bounded.

The second and third terms can be written as scaled sums of elementary operators,\\ $\big\{\lp E_{(j_1,0)},  2^{-(j_1,0)}\rp ;  j_1>0\big\}$ and $\big\{\lp E_{(0,j_2)}, 2^{-(0,j_2)}\rp ; j_2>0\big\}$ respectively, using the methods in the proof of \eqref{eq [T,p] decomp} so that:
\begin{align*}
    \sum_{\substack{ j_1> 0 \geq j_2}}  \ \eta [D_j,\psi_1] \eta' = \sum_{j_1>0}  \Big( \sum_{\substack{0\geq j_2}} \eta [D_j,\psi_1] \eta' \Big)  &=:\sum_{j_1>0} E_{(j_1,0)}\\
    \sum_{\substack{j_1 = j_2 > 0}} \eta  [D_j,\psi_1] \eta' = \sum_{j_2>0} 2^{-j_2} \Big( 2^{j_2} \eta [D_{(j_2,j_2)}, \psi_1] \eta' \Big) &=: \sum_{j_2 >0} 2^{-j_2}  E_{(0,j_2)}.
\end{align*}
Combining the results above, we conclude the proof of the claim and conclude the proof of Lemma \ref{lemma base case FK}. 
\end{proof}

\begin{proof}[Proof of Theorem \ref{main thm PK} for flag kernels]
Flag kernels on a direct product space form a subalgebra of product kernels. By Theorem \ref{main thm PK}, our new inverse kernel $L$ is thus a product kernel. To show that $L$ is in fact a flag kernel, we make a reduction using appropriate dilations. 

The flag kernel seminorms remain unchanged when we conjugate the operator $T$ with dilations $D_{(R_1, R_2)}TD_{(R_1^{-1}, R_2^{-1})}$, provided $R_1 \geq R_2$. In addition, if $|t_2|_2 \geq |t_1|_1$, then the growth conditions for product and flag kernels are equivalent. We thus further reduce to proving the growth condition in the case where $|t_2|_2 \leq |t_1|_1$. By homogeneity considerations, it suffices to show that
\begin{equation}\label{eq reduction FK}
    \sup_{|t_1|_1 \sim 1; \ |t_1|_1 \geq |t_2|_2} |\p^{\ap}L(t)| \lesssim 1.
\end{equation}
By conjugating the inverse operator with dilations $D_{(R_1, R_1)}$, the result follows. By retracing the proof of the cancellation condition for product kernels, we obtain the flag kernels cancellation condition for $L$ after a few straightforward modifications. Thus concluding the proof of Theorem \ref{main thm PK}. 
\end{proof}

\vspace{.2in}

\subsection*{Acknowledgements}
I would like to thank my advisor Brian Street for proposing this problem and for the many illuminating discussions that guided my investigation. This material is based upon work supported by the NSF under Grant No. DMS-2037851.

\printbibliography[title={References}]

\end{document}